\documentclass[oneside,twocolumn]{article}

\usepackage{cite}

\RequirePackage{fix-cm}
\usepackage{layout}
\usepackage[pdftex]{graphicx}
\def\isLong{false}
\newcommand{\myvspace}[1]{\vspace{0.0in}}

\usepackage{xspace}
\usepackage{helvet}
\usepackage{courier}
\usepackage{fancyhdr}
\usepackage{type1cm}         
\usepackage{subeqnarray}
\usepackage{graphicx}        
\usepackage{enumerate}
\usepackage{amsmath, latexsym}
\usepackage{array}
\usepackage{tabularx}
\usepackage{url}
\usepackage{relsize}
\usepackage{ifthen}
\usepackage{xcolor}
\usepackage[newitem,newenum]{paralist}
 
\def\isBeamer{false}

\newcommand{\mychapterskip}{\ifthenelse{\equal{\isMPK}{true}}
{\vspace{-.8in}}{}}

\newcommand{\fullversioncolor}{\color{black!50!blue}}
\newcommand{\versions}[2]{\ifthenelse{\equal{\isfullversion}{true}}{{\fullversioncolor#1}}{#2}}

\renewcommand{\vec}[1]{\mathbf{#1}}
\newcommand{\quot}[1]{``#1''}

\newcommand{\R}[1]{\ensuremath{\mathbb{R}^{#1}\xspace}}
\newcommand{\Euc}[1]{\ensuremath{\mathbf{E}^{#1}\xspace}}
\newcommand{\Eucg}[1]{\ensuremath{E(#1)\xspace}}
\newcommand{\Eucgd}[1]{\ensuremath{E^+(#1)\xspace}}

\newcommand\proj[1]{\ensuremath{\mathbf{P}(#1)}\xspace}
\newcommand\RP[1]{\ensuremath{\mathbb{R}{P^{#1}}\xspace}}

\newcommand{\e}[1]{\vec{e}_{#1}}
\newcommand{\EE}[1]{\vec{E}_{#1}}
\newcommand{\one}{\vec{1}}
\newcommand{\eye}{\vec{I}}
\newcommand{\inert}{\mathbf{A}}

\newcommand{\MN}{metric-neutral\xspace}

\newcolumntype{Y}{X}

\newcommand{\bivo}[1]{{\vec{#1}}}
\newcommand{\pip}{\bivo{\Xi}}

\newcommand{\pvelo}{\bivo{\Gamma}}
\newcommand{\velo}{\bivo{\Omega}}
\newcommand{\momo}{\bivo{\Pi}}

\newcommand{\sigo}{\bivo{\Sigma}}

\newcommand{\grade}[2]{\langle #1 \rangle_{#2}}

\newcommand{\pclal}[3]{\proj{\mathbb{R}_{#1,#2,#3}}\xspace}
\newcommand{\pdclal}[3]{\proj{\mathbb{R}^*_{#1,#2,#3}}\xspace}

\newcommand{\pdclplus}[3]{\proj{\mathbb{R}^{*+}_{#1,#2,#3}}\xspace}

\newcommand{\myboldhead}[1]{\vspace{0in}\hspace{0in}\textbf{#1.}}

\newcommand{\app}[1]{\Sec{sec:J}}

\definecolor{mypurple}{rgb}{1,0,1}

\definecolor{mygreen}{rgb}{0, .5, 0}

\newcommand{\Fig}[1]{Fig.~\ref{#1}}
\newcommand{\Tab}[1]{Table~\ref{#1}}
\newcommand{\Sec}[1]{Sect.~\ref{#1}}

\ifthenelse{\equal{\isBeamer}{true}}
{

}
{

}

\newcommand{\mycorrection}[1]{}

\renewcommand{\myboldhead}[1]{\vspace{.1in}\hspace{-.0in}\textbf{#1.}}

\ifthenelse{\equal{\isLong}{true}}
{
\title{Projective geometric algebra: \\New ways for graphics, games, and geometry}
\author{Charles G. Gunn\\Raum+Gegenraum\\Falkensee, Germany\\projgeom@gmail.com}
}
{
\title{Projective geometric algebra: \\ A new framework for doing euclidean geometry}
\author{Charles G. Gunn\\Raum+Gegenraum\\Falkensee, Germany\\projgeom@gmail.com}
}
\date{ }

\renewcommand{\myboldhead}[1]{\vspace{.1in}\hspace{-.14in}\textbf{#1.}}

\usepackage{layout}
\usepackage{hyperref}
\usepackage{paracol}

\graphicspath{{.}{./pictures/} }  

\usepackage[english]{babel}	
\usepackage{amsfonts}
\usepackage{latexsym}
\usepackage[protrusion=true,expansion=true]{microtype}				

\date{}

\begin{document}
\maketitle

%
\begin{abstract}

\ifthenelse{\equal{\isLong}{true}}
{
Projective geometric algebra (PGA) is a powerful, accessible toolkit for doing euclidean geometry in dimensions 2 and 3. It first establishes a powerful feature set that sets PGA apart from the standard approach using vector and linear algebra and analytic geometry, before showing how PGA can be profitably applied to current challenges in computer graphics,  computer vision, and related fields.
}
{
We introduce \emph{projective geometric algebra} (PGA), a modern, coordinate-free framework for doing euclidean geometry  featuring: uniform representation of points, lines, and planes; robust, \quot{parallel-safe} join and meet operations; compact, polymorphic syntax for euclidean formulas and constructions; a single intuitive \quot{sandwich} form for isometries; native support for automatic differentiation; and tight integration of kinematics and rigid body mechanics. Inclusion of vector, quaternion, dual quaternion, and exterior algebras as sub-algebras simplifies the learning curve and transition path for experienced practitioners. On the practical side, it can be efficiently implemented, while its rich syntax enhances programming productivity. The basic ideas are introduced in the 2D context; the 3D treatment focus on selected topics. Advantages to traditional approaches are collected in a table at the end. The article aims to be a self-contained introduction for practitioners of euclidean geometry and includes numerous examples, figures, and tables.
}
\end{abstract}


\section{Problem statement} What is the best representation  for doing euclidean geometry on computers? 
This question is a fundamental one for practitioners in an ever-growing range of application areas including computer graphics, computer vision, 3D games, virtual reality, robotics, CAD, animation, geometric processing, and discrete geometry. 
While available programming languages change and develop with reassuring regularity, the underlying geometric representations   
tend to be based on \textbf{v}ector and \textbf{l}inear \textbf{a}lgebra and \textbf{a}nalytic \textbf{g}eometry  (VLAAG for short), a framework that has remained virtually unchanged for 100 years.  
The article introduces \emph{projective geometric algebra} as a modern alternative for doing euclidean geometry (\cite{gunnthesis}, \cite{gunn2011}, \cite{gunn2017a}, \cite{gunn2017b}) and establishes that it enjoys significant advantages over VLAAG both conceptually and practically. 



\section{Feature list for doing euclidean geometry}
\label{sec:wishlist}
The standard  approach (VLAAG) has proved itself to be a robust and resilient toolkit. Countless engineers and developers use it to do their jobs. Why should they look elsewhere for their needs?  On the other hand, long-time acquaintance and habit can blind craftsmen to limitations in their tools, and subtly restrict the solutions that they look for and find. Many programmers have had an \quot{aha} moment when learning how to use the quaternion product to represent rotations without the use of matrices, a representation in which the axis and strength of the rotation can be directly read off of the four quaternion coordinates rather than laboriously extracted from the 9 entries of the matrix.

In the spirit of such \quot{aha!} moments we propose here a feature list for doing euclidean geometry on the computer that represents a significant advance over the features that VLAAG offers. We believe all developers will benefit from a framework that:
\begin{compactitem}
\item  is \textbf{coordinate-free}, 
\item  has a \textbf{uniform representation for points, lines, and planes}, 
\item  can calculate \textbf{meet and join} of  these geometric entities, while handling parallel elements correctly, 
\item  provides \textbf{compact  expressions} for  all classical euclidean formulas and constructions, including distances and angles, perpendiculars and parallels, orthogonal projections, and other metric operations, 
\item  has a \textbf{single, geometrically intuitive form} for euclidean motions, 
\item provides \textbf{automatics differentiation} of functions of one or several variables,
\item  provides a compact, efficient \textbf{model for kinematics and rigid body mechanics},  
\item lends itself to \textbf{efficient, practical implementation}, and
\item is \textbf{backwards-compatible} with existing representations such as vector, quaternion, dual quaternion, and exterior algebras.
\end{compactitem}

\subsection{Structure of the article}
In the rest of the article we will introduce geometric algebra in general and PGA in particular, on the way to showing that PGA in fact fulfills the above feature list.
Our treatment will be devoted to dimensions $n=2$ and $n=3$, the cases of most practical interest, and focuses on examples; readers interested in theoretical foundations are referred to the bibliography.  \Sec{sec:wiga} introduces geometric algebra and the associated \emph{geometric product} briefly before presenting three worked-out examples of PGA in action. \Sec{sec:roots} presents the historical background necessary to understand PGA. \Sec{sec:eucplane} then turns to PGA for the euclidean plane, written $\pdclal{2}{0}{1}$ where it introduces many of the fundamental features of PGA in this simplified setting: products of pairs of elements, formula factories from associativity, representation of isometries using sandwiches, and automatic differentiation.  \Sec{sec:eucspace} introduces PGA for euclidean 3-space; space restrictions limit this to a sketch of the role of bivectors, culminating in the Euler equations for rigid body motion expressed in the geometric algebra.  The rest of the article discusses implementation issues and compares the results with alternative approaches, notably VLAAG.

\myboldhead{What \emph{euclidean} is and isn't} First, we clarify what we mean by \quot{euclidean} since experience has shown this can be a stumbling block to approaching PGA. When we say \emph{doing euclidean geometry} we are referring to the geometry of \emph{euclidean space }$\Euc{n}$, not the \emph{euclidean vector space} $\R{n}$.  The elements of $\Euc{n}$ are points, those of $\R{n}$ are vectors; the motions of $\Euc{n}$ include translations and rotations, those of $\R{n}$ are rotations preserving the origin $O$. $\Euc{n}$ is intrinsically more complex than $\R{n}$: it is a differentiable metric space whose tangent space at each point is $\R{n}$.  We will see that euclidean PGA includes both $\Euc{n}$ and $\R{n}$ in an organic whole.

\section{What is geometric algebra?}
\label{sec:wiga}

PGA and VLAAG share common roots in classical 19th century mathematics which we describe in more detail in Sec. \ref{sec:roots}.  
The main idea behind geometric algebra is that \textbf{geometric primitives behave like numbers} -- they can be added and multiplied, have inverses and  appear in algebraic equations and functions.    The resulting interplay of algebraic and geometric aspects produces a synergy that has begun to attract the attention of applied mathematicians, see for example the textbook \cite{dfm07}.  PGA is a relative newcomer to the applied geometric algebra scene: the idea appeared in the modern literature first in \cite{selig00} and was given its name in \cite{gunn2017a}. 

Fortunately, many features of PGA are already familiar to practitioners.  It is based on homogeneous coordinates, widely used in computer graphics, and it contains within it classical vector algebra, as well as the quaternion and dual quaternion algebras, increasingly popular tools for modeling kinematics and mechanics.  The exterior algebra, a powerful structure that models the subspaces of $\R{n}$ (or projective space $\RP{n}$), is also contained as a sub-algebra.  
PGA in fact can be compared to a whole organism in which each of these sub-algebras first finds its true place in the scheme of things.  For readers familiar with the use of conformal geometric algebra (CGA), a detailed comparison of PGA and CGA for doing euclidean geometry may be found in \cite{gunn2017a} where the two algebras are assigned complementary positions in the GA eco-system.  The focus of our comparison here is with VLAAG restricted to flat primitives (e. g., points, lines, and planes).  Before turning to the formal details we present three examples of  PGA at work, solving tasks in 3D euclidean geometry, to give a flavor of actual usage. Readers who prefer a more systematic introduction are encouraged to skip over to \Sec{sec:gralg}.

\subsection{\textbf{Example 1}: Working with lines and points in 3D}
\begin{quote}
\textbf{Task:} Given a point $\vec{P}$ and a non-incident line $\momo$ in $\Euc{3}$, find the unique line $\sigo$ passing through $\vec{P}$ which meets $\momo$ orthogonally.
\end{quote}
  \begin{figure}[h]
  \begin{centering}
 \def\xyz{.48}
 \def\zyx{.01}
{\setlength\fboxsep{0pt}\fbox{\includegraphics[width=\xyz\columnwidth]{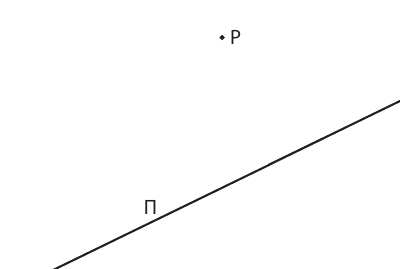}}}\hspace{\zyx\columnwidth}
{\setlength\fboxsep{0pt}\fbox{\includegraphics[width=\xyz\columnwidth]{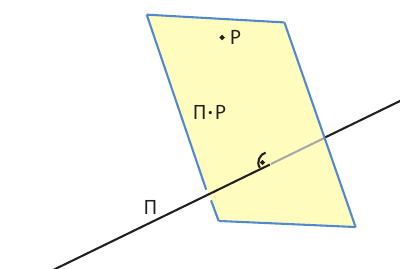}}}\\ \vspace{\zyx\columnwidth}
{\setlength\fboxsep{0pt}\fbox{\includegraphics[width=\xyz\columnwidth]{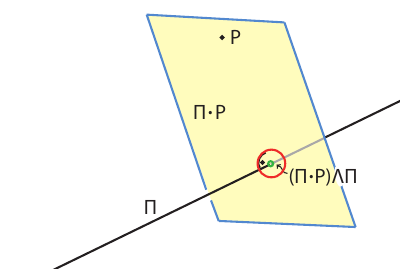}}}\hspace{\zyx\columnwidth}
{\setlength\fboxsep{0pt}\fbox{\includegraphics[width=\xyz\columnwidth]{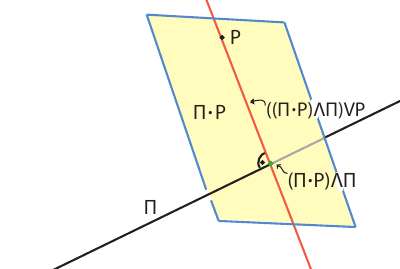}}}
\caption{Geometric construction in PGA.} 
\label{fig:perpPtLn}
\end{centering}
\end{figure}
\myvspace{-.1in}
In PGA,  geometric primitives such as points, lines, and planes, are represented by vectors of different \emph{grades}, just as in an exterior algebra.  A plane is a 1-vector,  a line is a 2-vector, and a point is a 3-vector. (A scalar is a 0-vector; we'll meet 4-vectors in \Sec{sec:3dscrew}).  Hence the algebra is called a \emph{graded} algebra. The geometric relationships between primitives is expressed via  the \emph{geometric product}, an associative bilinear product defined on these $k$-vectors. 

The geometric product $\vec{P} \momo$, for example, of a point $\vec{P}$ (a 3-vector) and a line $\momo$ (a 2-vector) consists two pieces: the plane perpendicular to $\momo$ passing through $\vec{P}$ (a 1-vector, written as $\momo \cdot  \vec{P}$), and the normal direction to the plane spanned by $\vec{P}$ and $\momo$ (a 3-vector, written $\vec{P} \times \momo$).  The sought-for line  $\sigo$ can then be constructed as shown in \Fig{fig:perpPtLn}: 
 \begin{compactenum}
 \item $\momo \cdot \vec{P}$ is the plane through $\vec{P}$ perpendicular to $\momo$,
 \item $(\momo \cdot \vec{P} )\wedge \momo)$ is the meet ($\wedge$) of $\momo \cdot \vec{P}$ with $\momo$,
 \item $\sigo := ((\momo \cdot \vec{P} )\wedge \momo)  \vee \vec{P}$ is the join ($\vee$) of this point with $\vec{P}$.
 \end{compactenum}
 The meet ($\wedge$) and joint ($\vee$) operators are part of the geometric algebra and are discussed in more detail below in \Sec{sec:gralg}.
 
   \begin{figure}[h]
   \begin{centering}
 \def\xyz{.34}
 \def\xyw{.68}
 \def\zyx{.01}
{\setlength\fboxsep{0pt}\fbox{\includegraphics[height=\xyz\columnwidth]{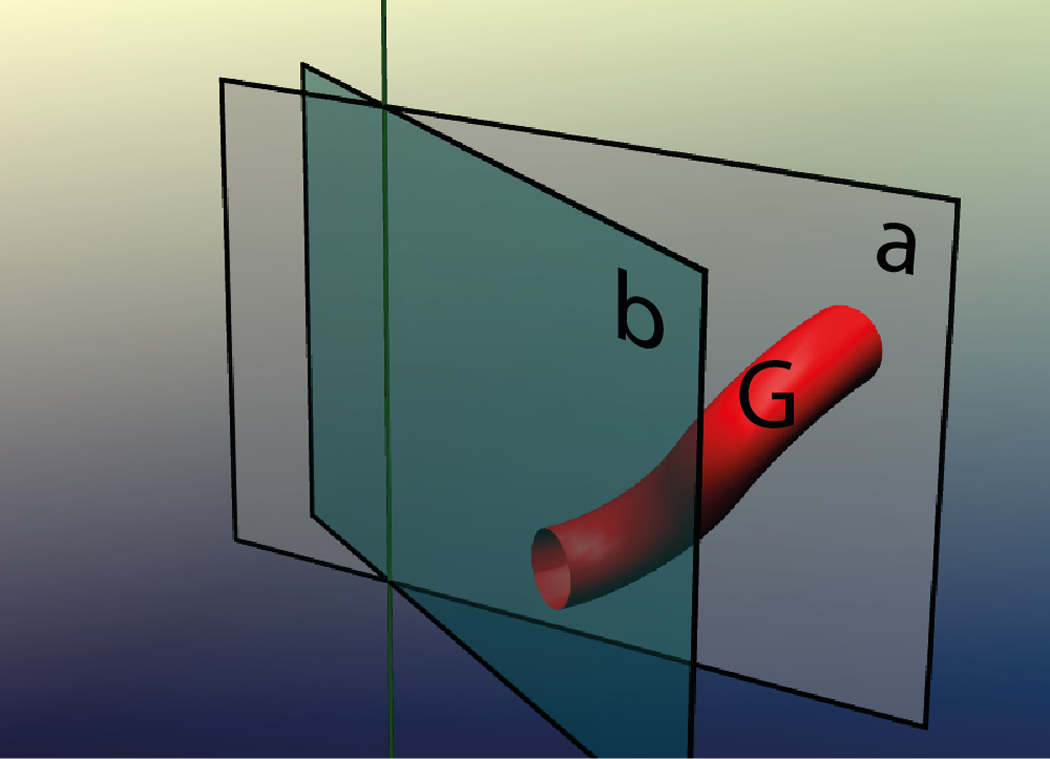}}}\hspace{\zyx\columnwidth}
{\setlength\fboxsep{0pt}\fbox{\includegraphics[height=\xyz\columnwidth]{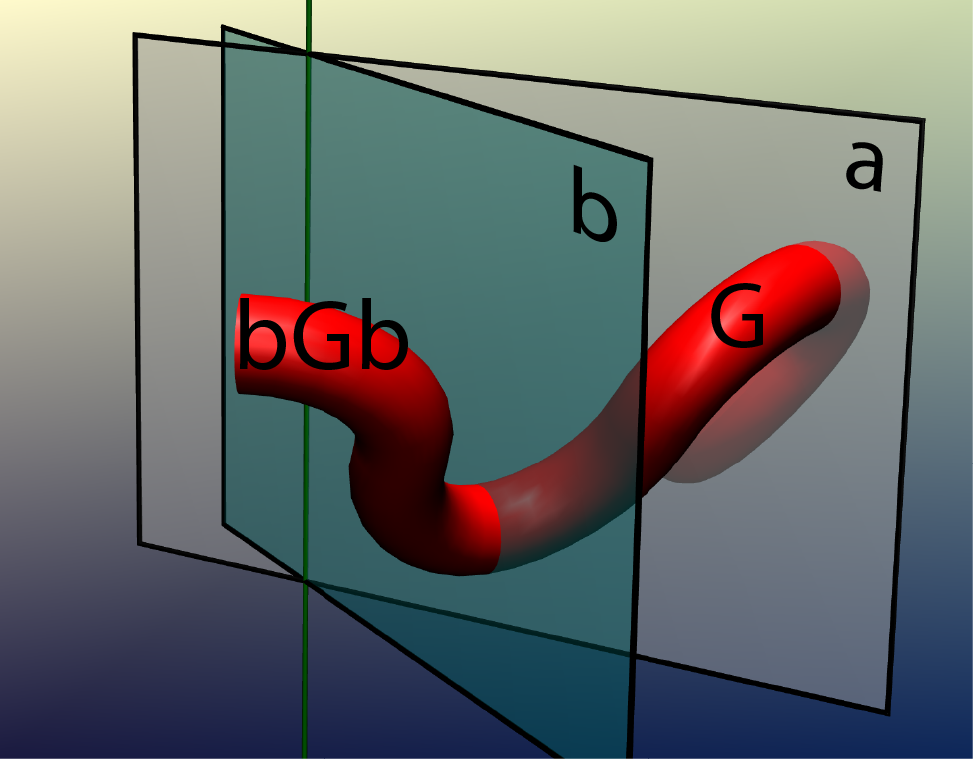}}} \\ \vspace{\zyx\columnwidth}
{\setlength\fboxsep{0pt}\fbox{\includegraphics[height=\xyw\columnwidth]{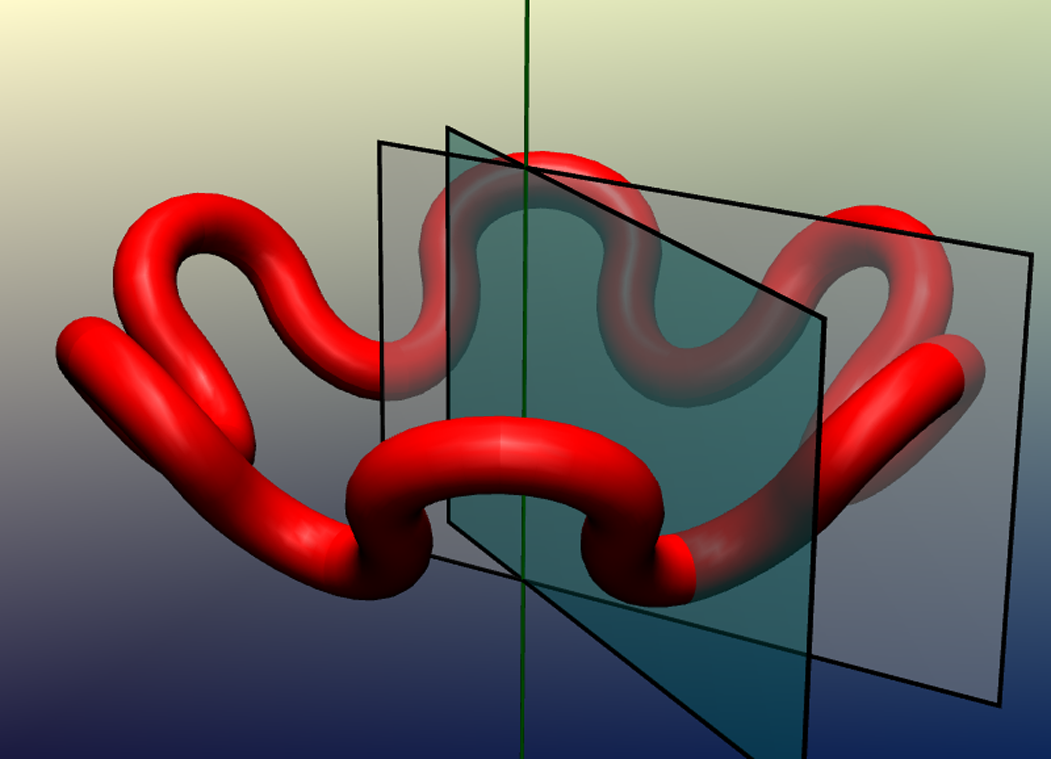}}}
\caption{Creating a 3D kaleidoscope in PGA using sandwich operators.} 
\label{fig:3dkal}
\end{centering}
\end{figure}
\myvspace{-.4in}

\subsection{\textbf{Example 2}: A 3D Kaleidoscope}
\label{sec:3dkal}
\begin{quote}
\textbf{Task:} A $k$-kaleidoscope is a pair of mirror planes $\vec{a}$ and $\vec{b}$  in $\Euc{3}$ that meet at an angle $\frac{\pi}{k}$. Given some geometry $\vec{G}$ generate the view of $\vec{G}$ seen in the kaleidoscope.
\end{quote}

In PGA, the reflection in a plane $\vec{a}$ (a 1-vector) is represented using the geometric product by the \quot{sandwich} operator $\vec{a}\vec{G}\vec{a}$ (where $\vec{G}$ may be any $k$-vector). See  \Fig{fig:3dkal}. The left-most image shows the given situation, where $\vec{G}$ is a red tube (modeled by any combination of 1-, 2-, and 3-vectors) stretching between the two planes.  The middle image shows the result of applying the sandwich $\vec{b}\vec{G}\vec{b}$ to the geometry (behind plane $\vec{a}$ one can also see $\vec{a}\vec{G}\vec{a}$, unlabeled).  Note that we can and do normalize the plane $\vec{a}$ to satisfy $\vec{a}^2=1$ (where $\vec{a}^2$ is  the geometric product of $\vec{a}$ with itself), which is consistent with the fact that repeating a reflection yields the identity.  The right image shows the result of applying all possible alternating products of these two reflections to $\vec{G}$ (e. g., $\vec{b}\vec{a}\vec{G}\vec{a}\vec{b}$, etc.). Since the mirrors meet at the angle $\frac{\pi}{6}$, this process closes up in a ring consisting of 12 copies of $\vec{G}$. (To be precise, $(\vec{a}\vec{b})^6 = (\vec{b}\vec{a})^6 = 1$). 
\myvspace{-.1in}
\subsection{\textbf{Example 3}: A continuous 3D screw motion}
\label{sec:3dscrew}
\begin{quote}
\textbf{Task:}  Represent a continuous screw motion in 3D.
\end{quote}
Recall that the general orientation-preserving isometry of $\Euc{3}$ is a \emph{screw motion}, which rotates around a unique fixed line (the  \emph{axis}) while translating parallel to it.  The ratio of the translation distance to the angle of rotation (in radians) is called the \emph{pitch} of the screw motion.  A rotation has pitch 0, and translation has pitch \quot{$\infty$}.

  \begin{figure}[t]
 \def\xyz{.34}
 \def\xyw{.69}
 \def\zyx{.01}
 \begin{centering}
{\setlength\fboxsep{0pt}\fbox{\includegraphics[height=\xyz\columnwidth]{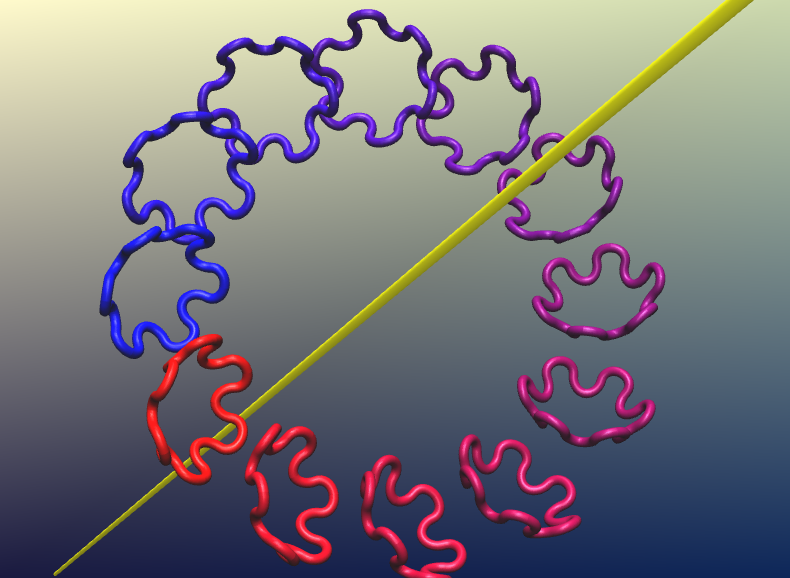}}}\hspace{\zyx\columnwidth}
{\setlength\fboxsep{0pt}\fbox{\includegraphics[height=\xyz\columnwidth]{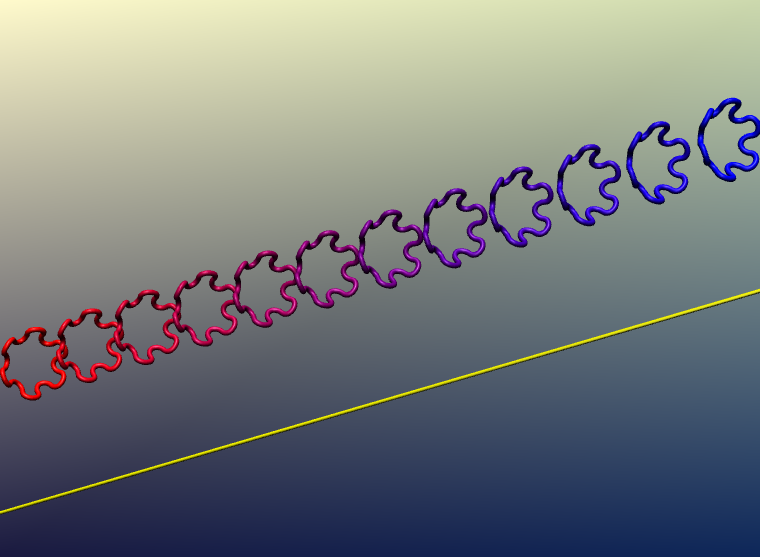}}} \\ \vspace{\zyx\columnwidth}
{\setlength\fboxsep{0pt}\fbox{\includegraphics[height=\xyw\columnwidth]{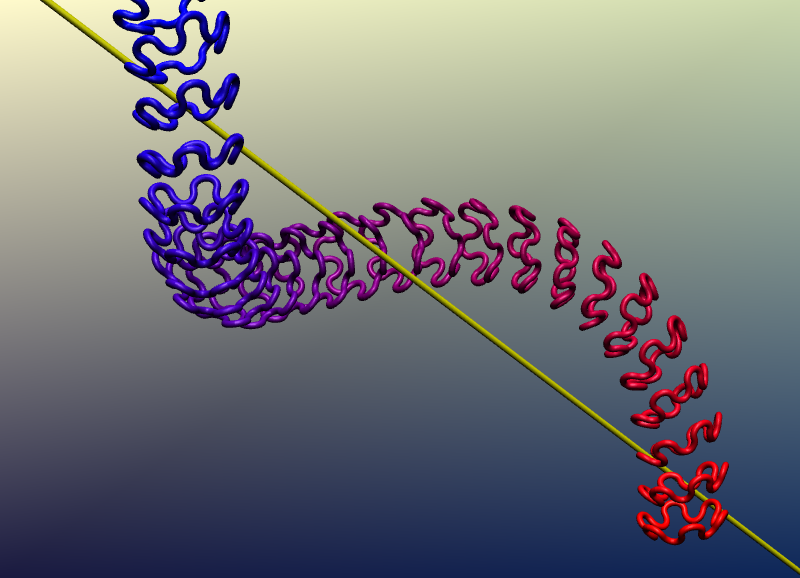}}}
\caption{Continuous rotation, translation, and screw motion in PGA by exponentiating a bivector.} 
\label{fig:screwMotion}
\end{centering}
\end{figure}

We first show how to represent a rotation in PGA.  A line in $\Euc{3}$, passing through the point $\vec{P}$ with direction vector $\vec{V}$, is given by the join operation $\velo = \vec{P} \vee \vec{V}$ (yellow line in Fig. \ref{fig:screwMotion}). To obtain the  \emph{rotation} around $\velo$ of angle $\alpha$  define the \emph{rotor} $e^{t \velo}$.  The exponential function is  evaluated using the geometric product in the formal power series of $e(x)$; it behaves like the imaginary exponential $e^{ti}$ since  we can and do normalize $\velo$  to satisfy $\velo^2=-1$. Then the continuous rotation around $\velo$ applied to $\vec{G}$ is given by the \emph{sandwich operator} $e^{t \velo}\vec{G} e^{-t \velo}$.  At $t=0$ it is the identity; and at $t=\frac{\alpha}{2}$ it represents the rotation of angle $\alpha$ around $\velo$. See the left image above, which shows the result for a sequence of $t$-values between $0$ and $\pi$. Readers familiar with the quaternion representation of rotations should recognize these formulas.

To obtain instead a translation, we apply the \emph{polarity} operator of PGA to $\velo$ to produce $\velo^\perp$, the \emph{orthogonal complement} of $\velo$.  $\velo^\perp$ is an \emph{ideal} line, or so-called \quot{line at infinity}, in this case consisting of the directions of all lines which meet $\velo$ at right angles.  
It is obtained from $\velo$ by multiplying by a special 4-vector, the unit \emph{pseudoscalar} $\eye$: $\velo^\perp := \velo \eye$.  A continuous translation parallel to $\velo$ is then given by a sandwich with the \emph{translator} $e^{t \velo^\perp} $.  See the middle image above. 

Let the pitch of the screw motion be $p\in \mathbf{R}$.  Then the desired screw motion is given by a sandwich operator with  the \emph{motor}  $e^{t(\velo+p\velo^\perp)}$.  Note that the exponent is a linear combination of the (commuting) rotational and translational exponents, weighted by $p$. 
See image on the right above.

We hope these examples have whetted your appetite to explore further. We now turn to a quick exposition of the history of PGA, introducing the essential ideas needed to master the subject.
\myvspace{-.1in}
\section{Roots of PGA}
\label{sec:roots}
\begin{figure}[b]
 \def\xyz{.3}
 \def\zyx{.02}
{\setlength\fboxsep{0pt}{\includegraphics[height=\xyz\columnwidth]{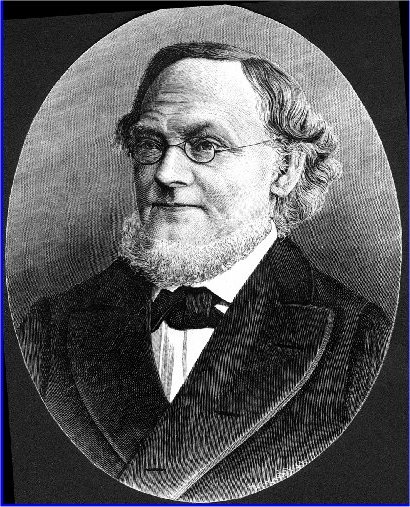}}}\hspace{\zyx\columnwidth}
{\setlength\fboxsep{0pt}{\includegraphics[height=\xyz\columnwidth]{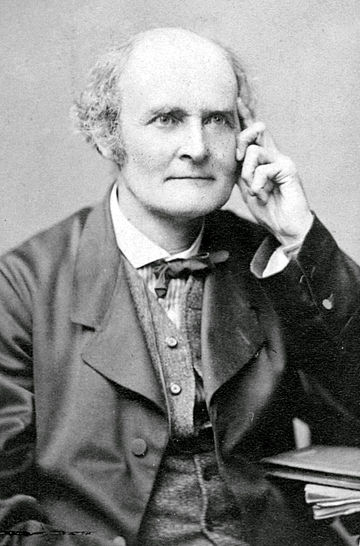}}}\hspace{\zyx\columnwidth}
{\setlength\fboxsep{0pt}{\includegraphics[height=\xyz\columnwidth]{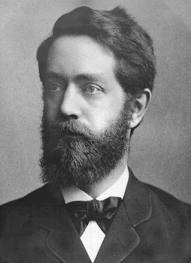}}}\hspace{\zyx\columnwidth}
{\setlength\fboxsep{0pt}{\includegraphics[height=\xyz\columnwidth]{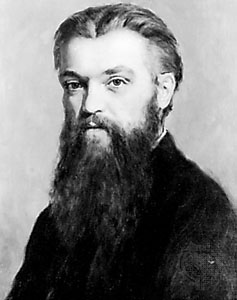}}}
\caption{Important figures in the development of PGA (l. to r.): Hermann Grassmann (1809-1877), Arthur Cayley (1821-1895), Felix Klein (1849-1925), William Clifford (1845-1879).} 
\label{fig:inventors}
\end{figure}

Both the standard approach to doing euclidean geometry and the geometric algebra approach described here can be traced back to 16th century France. 
The analytic geometry of Ren\'{e} Descartes (1596-1650) leads to the standard toolkit used today.   His contemporary and friend Girard Desargues (1591-1661), an architect, confronted with the riddles of the newly-discovered perspective painting, invented \emph{projective geometry}, containing additional, so-called \emph{ideal}, points where parallel lines meet.  Projective geometry is characterized by a deep symmetry called \emph{duality}, that asserts that every statement in projective geometry has a dual partner statement, in which, for example, the roles of point and plane, and of join and intersect, are exchanged. More importantly, the truth content of a statement is preserved under duality. We will see below that duality plays an important role in PGA. 

Mathematicians in the 19th century (Cayley and Klein) showed how, using an algebraic structure called a \emph{quadratic form},  the euclidean metric could be built back into projective space.  It is this \emph{Cayley-Klein model }of euclidean geometry that forms the backbone of PGA. 
We focus to start with on the euclidean plane, where we can get to know all the interesting behavior before moving on to $n=3$, where the real interest lies, but also much more complex behavior.    We next turn to how to represent the subspaces of projective space in an algebraic structure.

\subsection{Exterior algebra of subspaces of $\RP{n}$}
\label{sec:gralg}
The backbone of every geometric algebra is an \emph{exterior} algebra.  Exterior algebras were, like so many other results in this field, discovered by Hermann Grassmann (\cite{grassmann44}) and are sometimes called \emph{Grassmann} algebras. An exterior algebra (as used here) mirrors the subspace structure of projective space $\RP{n}$. One can build up the subspaces of  projective $n$-space  by \emph{joining points}; or  by \emph{intersecting hyperplanes}.   Duality ensures that these two approaches are completely equivalent and neither \emph{a priori} is to be preferred.  Each construction produces a separate exterior algebra. 

 In a {standard} exterior algebra G, the elements of grade $k$ for $k = 1, 2, ... n$,  represent the subspaces of dimension $k-1$.  For example, for $n=2$, the 1-vectors are points, and the 2-vectors are lines.    The graded algebra also has  elements of grade 0, the  scalars (the real numbers $\mathbb{R}$); and elements of grade $(n+1)$ (the highest non-zero grade), the \emph{pseudoscalars}. 
  All elements of the exterior algebra have projective coordinates; so each has a non-zero weight which can be freely chosen and, as we will see, often expresses important geometric information. The space of $k$-vectors is a vector space written $\bigwedge^{k}$.

  \begin{figure}[h]
   \centering
 \def\xyz{.45}
 \def\zyx{.05}
{\setlength\fboxsep{0pt}{\includegraphics[height=\xyz\columnwidth]{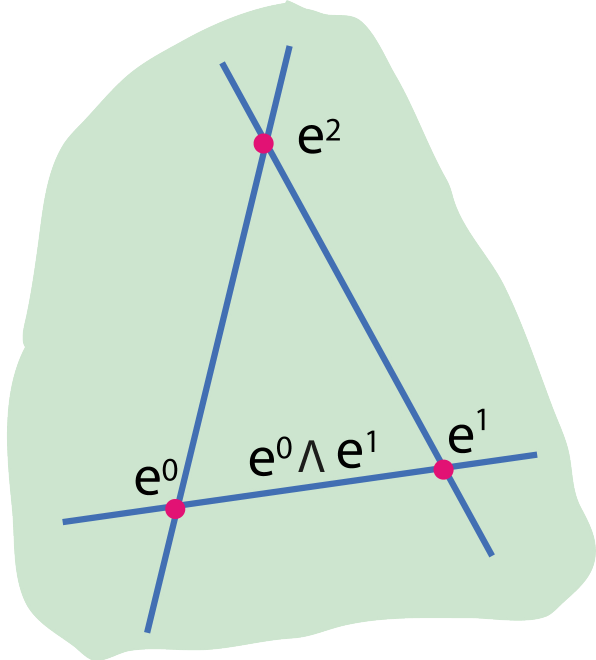}}}\hspace{\zyx\columnwidth}
{\setlength\fboxsep{0pt}{\includegraphics[height=\xyz\columnwidth]{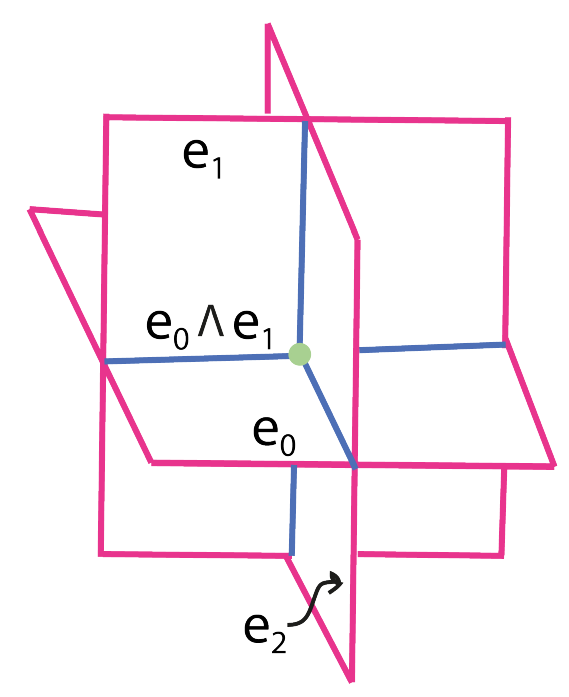}}}
\caption{\emph{Left:} The plane $\vec{e}^0\wedge\vec{e}^1\wedge\vec{e}^2$ (green) created by joining 3 points in G, the standard exterior algebra (written with raised indices).  \emph{Right:} The meeting point $\vec{e}_0\wedge\vec{e}_1\wedge\vec{e}_2$ (green) of three planes in G*, the dual exterior algebra (written with lowered indices).} 
\label{fig:dualAufbau}
\end{figure}
\myvspace{-.1in}
  There is an anti-symmetric associative bilinear  product, the \emph{outer}, or \emph{wedge}, {product}, in an exterior algebra. It represents, in the standard exterior algebra,  the subspace obtained by \emph{joining} the two arguments.  The outer product of a linearly independent $k$- and $m$-vector is the $(k+m-1)$-dimensional subspace that they span, otherwise it is 0.   In the \emph{dual} exterior algebra G*, on the other hand,  elements of grade $k$ represent the subspaces of dimension $n-k$, and the outer product is the \emph{meet} {operator}.  
 Consult \Fig{fig:dualAufbau}, which shows how the wedge product of three points in G is a plane, while the wedge product of three  planes in G* is a point.
 Using the fact that every geometric entity occurs once in each exterior algebra, it's possible to \quot{import} the outer product from one algebra into the the other using a grade-reversing isomorphism called the \emph{Poincar\'{e} duality} map, so join and meet are available within a single algebra, see  \cite{gunnthesis}, \S 2.3.1.
 
In the following we focus on G*, since as we'll see below in \Sec{sec:ipeg}, euclidean PGA is built using G*, not G. 
We write the  outer product of G*, the meet operator, as $\wedge$, and the  join operator, imported from G, as $\vee$.  That's easy to remember due to their similarity to the set operations $\cap$ and $\cup$.  Important to note:
Working in projective space guarantees that the meet of  parallel lines and planes, as well as the join of euclidean and ideal elements, are handled seamlessly, without \quot{special casing} -- one of the features on our initial wish-list. Details lie outside the scope of this treatment.

\subsection{Adding an inner product }

The exterior algebra with its outer product(s) is a powerful tool for calculating parallel-safe incidence relationships.  But to calculate euclidean angles and distances, one needs additionally to introduce an \emph{inner product.}  
An inner product of dimension $n$ is a symmetric bilinear product on the space of 1-vectors, characterized by three non-negative integers $(p,m,z)$, called its \emph{signature}, with $p+m+z=n$, such that there is a basis in which the squares of the $n$ basis elements consist of $p$ +1's, $m$ -1's, and $z$ 0's.  The familiar positive definite inner product of the euclidean \emph{vector} space $\R{n}$ has the signature $(n,0,0)$. We'll discover the proper signature for euclidean space below in \Sec{sec:ipeg}. 

\myboldhead{Geometric product} We combine the  outer product ($\wedge$)  with an inner product ($\cdot$) to obtain a \textbf{geometric product} on the 1-vectors of an exterior algebra.  For two one-vectors $\vec{a}$ and $\vec{b}$ the geometric product takes the simple form: 
\begin{align*}
{\vec{a} \vec{b} := \vec{a} \cdot \vec{b} + \vec{a} \wedge \vec{b}}
\end{align*}
where the two terms on the right-hand side are a scalar and 2-vector, resp.
This geometric product on 1-vectors can be naturally extended to all $k$-vectors to produce an associative product defined on the whole exterior algebra (see \cite{lang71}). The algebra equipped with this geometric product is called a \emph{geometric algebra}. This is the name Clifford gave it when he introduced it in \cite{clifford78}. We use this term as a synonym for \emph{Clifford algebra}. 
Because we work in projective space we call it a \emph{projective} geometric algebra or PGA for short.  It uses $(n+1)$-dimensional coordinates to model euclidean geometry. This distinguishes it from VGA (vector geometric algebra), that is build on $n$-dimensional vector space coordinates; and CGA (conformal geometric algebra) which uses $(n+2)$-dimensional coordinates to model $n$-dimensional euclidean space (introduced in \cite{hlr01}, see also \cite{dfm07}).  
There are also non-euclidean versions of PGA for hyperbolic and elliptic space; interested readers can consult \cite{gunnthesis}.

\myboldhead{GA Terminology} In general, the product of a $k$-vector and an $m$-vector is, just like the product of two 1-vectors shown above,  a sum of components of different grades, each expressing a different geometric aspect of the product.  Such a general element is  called a \emph{multi-vector}. A multi-vector $\vec{M}$ can be written then as a sum of different grades: $\vec{M} = \sum_{i=0}^n\grade{\vec{M}}{i}$.  We can also write the above geometric product as:  $\vec{a} \vec{b} := \grade{\vec{a} \vec{b}}{0} +\grade{ \vec{a}  \vec{b}}{2}$. We define the lowest-grade part of the geometric product of a $k$-vector and an $m$-vector to be the \emph{inner product} and written $ \vec{a} \cdot \vec{b}$ (even when it's not a scalar); the $(k+m)$-grade part  coincides with $\wedge$.  In what follows, we introduce notation for other grade parts of the geometric product as the situation requires. A $k$-vector which can be written as the product of 1-vectors is called a \emph{simple} $k$-vector. We sometimes call 2-vectors \emph{bivectors}, and 3-vectors, \emph{trivectors}.

\subsection{An  inner product for euclidean geometry}
\label{sec:ipeg}
We now turn to determining the correct inner product for doing euclidean geometry.  This inner product, for example, reveals itself when calculating the angle between two lines in the plane. 
Let \[
a_{0}x + b_{0}y + c_{0} = 0,~~~~~~~
a_{1}x + b_{1}y + c_{1} = 0
\]  be two oriented lines which intersect at an angle $\alpha$.  We can assume without loss of generality that the coefficients satisfy ${a_{i}^{2} + b_{i}^{2} = 1}$. Then it is not difficult to show that $a_{0} a_{1} + b_{0} b_{1} = \cos{\alpha}$. 

\begin{figure}[h]
  \centering
  \def\xyz{0.8}
    \setlength\fboxsep{0pt}\fbox{\includegraphics[width=\xyz\columnwidth]{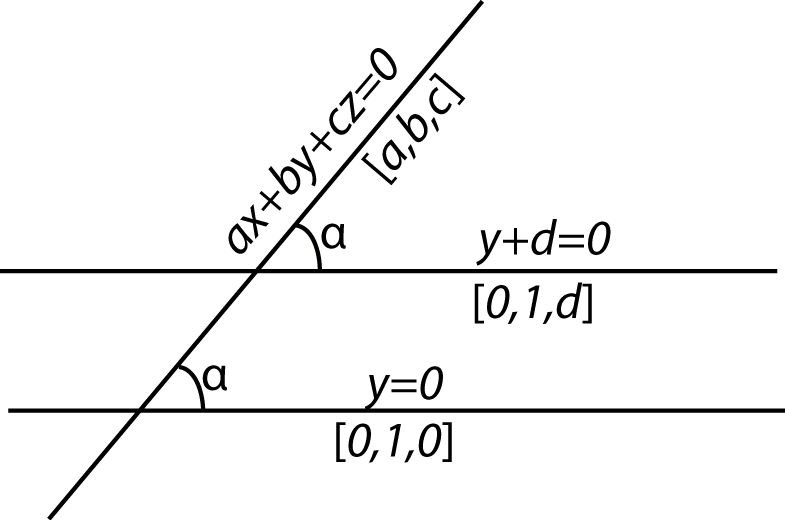}}
  \caption{Angles of euclidean lines. }
\label{fig:anglelines}
\end{figure} 

The third coordinate of the lines makes no difference in the angle calculation. Indeed,  translating a line changes only its third coordinate,  leaving  the angle between the lines unchanged.  Refer to \Fig{fig:anglelines} which shows an example involving a general line and a pair of horizontal lines. Hence the proper signature for measuring angles in $\Euc{2}$ is $(2,0,1)$.  
A similar argument applies in dimension $n$, yielding  the  signature  $(n,0,1)$ for $\Euc{n}$.  Such a signature, or metric, is called \emph{degenerate} since $z\neq 0$. The resulting geometric algebra is written $\pdclal{n}{0}{1}$.  The $*$ in the name says that the algebra is built on G*, the dual exterior algebra, since the inner product is defined on hyperplanes (lines in the case $n=2$) instead of points. 
$\pclal{n}{0}{1}$ models a qualitatively different metric space called \emph{dual euclidean space}.

PGA's development reflects the fact that much of the existing literature on geometric algebras deals only with non-degenerate metrics, and several misconceptions regarding degenerate metrics have become widespread. (See \cite{gunn2017a} for a thorough analysis and refutation of these misconceptions.)  After long experience we are convinced that the degenerate metric, far from being a liability, is the secret of PGA's success -- only a degenerate metric can model the metric relationships of euclidean geometry (see \cite{gunn2017a}, \S 5.3). 

%
\section{The euclidean plane $\pdclal{2}{0}{1}$}
\label{sec:eucplane}

We give now a brief introduction to PGA via euclidean plane geometry.  Readers eager to know more are referred to \cite{gunn2017b}. The approach presented here can be carried out in a coordinate-free way (\cite{gunn2017b}, Appendix).  But for an introduction it's easier and also helpful to refer occasionally to coordinates.  The coordinates we'll use are sometimes called \emph{affine coordinates} for euclidean space. We add an extra coordinate to standard $n$-dimensional coordinates, either a 1 (for euclidean points) or a 0 (for ideal points). 

\begin{figure}[h]
  \centering
  \def\xyz{0.95}
    \setlength\fboxsep{0pt}\fbox{\includegraphics[width=\xyz\columnwidth]{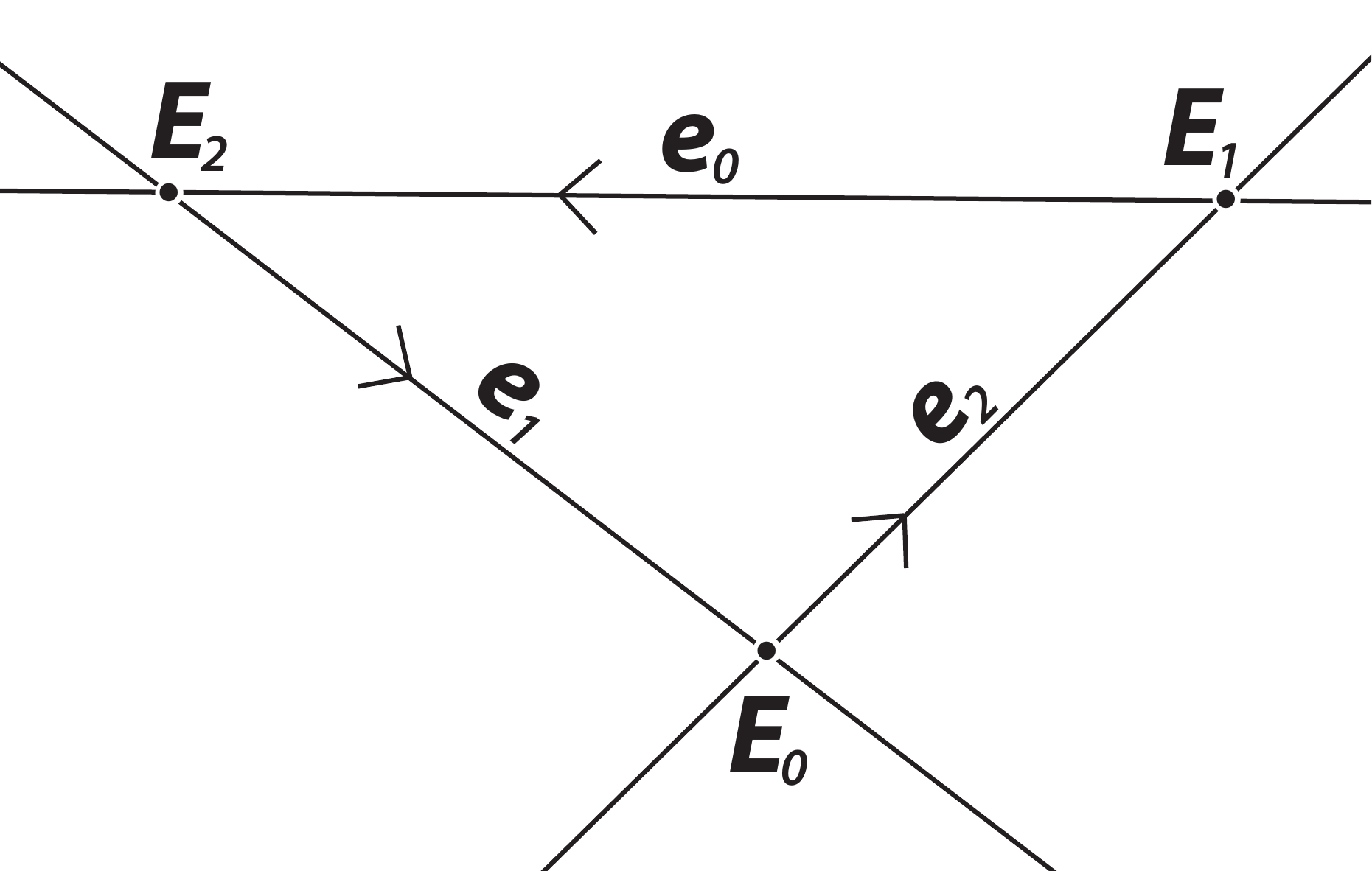}}
  \caption{Fundamental triangle of coordinate system. }
\label{fig:fundtri}
\end{figure} 

A perspective figure of the basis elements is shown in \Fig{fig:fundtri}.  The basis 1-vector $\e{0}$ represents the ideal line (sometimes called the \quot{line at infinity} and written $\omega$). $\e{1}$ and $\e{2}$ represent the coordinate lines $x=0$ and $y=0$, resp.  Note that for a 1-vector $\vec{a}$, $\vec{a}^2 = \vec{a}\cdot \vec{a}$ since $\vec{a}\wedge\vec{a}=0$ by anti-symmetry of $\wedge$. We choose then basis vectors to  satisfy $\e{0}^2 = 0$ and $\e{1}^2=\e{2}^2 = 1$. A basis for the 2-vectors is given by the intersection points of these orthogonal basis lines: \[\EE{0} := \e{1}\e{2},~~~ \EE{1}:=\e{2}\e{0}, ~~~\EE{2}:=\e{0}\e{1}\] whereby $\EE{0}$ is the origin, $\EE{1}$ and $\EE{2}$ are the $x$- and $y-$ directions (ideal points), resp.  They satisfy $\EE{0}^2=-1$ while $\EE{1}^2=\EE{2}^2=0$.  Finally, the unit pseudoscalar $\eye := \e{0}\e{1}\e{2}$ represents the whole plane and satisfies $\eye^2=0$.  The full 8x8 multiplication table of these basis elements can be found in \Tab{tab:cl201}.

\begin{table}[b]
\centering
\renewcommand{\arraystretch}{1.1}
\begin{tabularx}{\columnwidth} {| Y  || Y | Y  |  Y  | Y | Y | Y | Y | Y  |} \hline 
           & $\one$ & $\e{0}$ & $\e{1}$ & $\e{2}$ & $\EE{0}$ & $\EE{1}$ & $\EE{2}$ & $\eye$  \\ \hline \hline
$\one$        & $\one$ & $\e{0}$ & $\e{1}$ & $\e{2}$ & $\EE{0}$ & $\EE{1}$ & $\EE{2}$ & $\eye$  \\ \hline
$\e{0}$  & $\e{0}$ & $0$     & $\EE{2}$  & $-\EE{1}$ & $\eye$    & $0$      & $0$          & $0      $  \\ \hline
$\e{1}$  & $\e{1}$ & $-\EE{2}$ & $\one$& $\EE{0}$ & $\e{2}$ & $\eye$ & $-\e{0}$ & $\EE{1}$ \\ \hline
$\e{2}$  & $\e{2} $  & $\EE{1}$ & $-\EE{0}$ & $\one$ & $-\e{1}$ & $\e{0}$ & $\eye$ & $\EE{2}$ \\ \hline
$\EE{0}$  & $\EE{0}$ & $\eye$ & $-\e{2}$   & $\e{1}$   & $-\one$ & $-\EE{2}$ & $\EE{1}$ & $-\e{0}$  \\ \hline 
$\EE{1}$  & $\EE{1}$ & $0$     & $\eye$     & $-\e{0}$    & $\EE{2}$ & $0$ & $0$ & $0$  \\ \hline
$\EE{2}$  & $\EE{2}$ & $0$     & $\e{0}$   & $\eye$    & $-\EE{1}$ & $0$ & $0$ & $0$  \\ \hline
$\eye$    & $\eye$     & $0$     & $\EE{1}$ & $\EE{2}$ & $-\e{0}$ & $0$ & $0$ & $0$ \\ \hline
\end{tabularx}
\myvspace{.1in}
\caption{Multiplication table for the geometric product in $\pdclal{2}{0}{1}$}
\label{tab:cl201}
\end{table}

\subsection{Normalizing $k$-vectors}

Just as with euclidean vectors in $\R{n}$, it's possible and often preferable to normalize simple $k$-vectors.  For $n=2$, for example, we can normalize a euclidean line $\vec{a}$ so that $\vec{a}^2=1$; a euclidean point $\vec{P} = (x,y,1)$ is normalized and satisfies $\vec{P}^2 = -1$. This gives rise to a standard norm on euclidean $k$-vectors $\vec{X}$ that we write $\|\vec{X} \|$.  

\myboldhead{The ideal norm} Such a normalization is not possible for ideal elements, since these satisfy $\vec{X}^2=0$.  For ideal elements we define a second  norm $\| \|_\infty$, the \emph{ideal} norm.  In terms of the coordinates introduced above, for an ideal point $(x,y,0)$, $\|(x,y,0)\|_\infty = \sqrt{x^2+y^2}$. This agrees with the standard Euclidean vector space norm restricted to the  subspace satisfying $z=0$. (Note: A coordinate-free definition of the ideal norm of an ideal point $\vec{V}$ is given by $\|\vec{V}\|_\infty := \| \vec{V} \vee \vec{P}\|$ for any normalized euclidean point $\vec{P}$.) In fact, an ideal point is essentially a free vector, a fact already recognized by Clifford \cite{clifford73}.  We will see  that these two norms harmonize remarkably with each other, producing \emph{polymorphic} formulas -- formulas that produce correct results for any combination of euclidean and ideal arguments. We meet such an example in the product of two euclidean lines in the following section.

In the discussions below, we assume that all the arguments have been  normalized with the appropriate norm since, just as in $\R{n}$, it simplifies some discussions.  


\subsection{Examples: Products of pairs of elements in 2D}
\label{sec:prodpr}
\begin{figure}[h]
 \centering
 \def\xyz{.95}
{\setlength\fboxsep{0pt}\fbox{\includegraphics[width=\xyz\columnwidth]{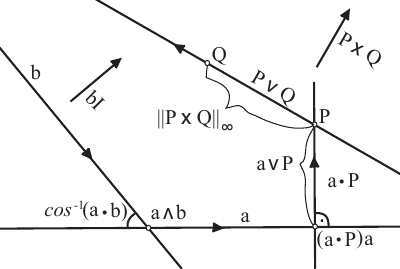}}}
\caption{Selected geometric products of pairs of simple vectors.}
\label{fig:geomprod}
\end{figure}

 We get to know the geometric product better by considering basic products.  
 A full discussion can be found in \cite{gunn2017b}.  Consult 
\Fig{fig:geomprod}.

\begin{compactitem}
\item \textbf{Multiplication by the pseudoscalar.} Multiplication by the pseudoscalar $\eye$ maps a $k$-vector onto its orthogonal complement with respect to the euclidean metric.  For a euclidean line $\vec{a}$, $\vec{a}^\perp := \vec{a}\eye$ is an ideal point perpendicular to the direction of $\vec{a}$.  For a euclidean point $\vec{P}$, $\vec{P} ^\perp := \vec{P}\eye$ is the ideal line $\vec{e}_0$. And $\eye^2 = 0$ is characteristic of degenerate metrics.

\item \textbf{Product of two euclidean lines.} 

$\vec{a} \cdot \vec{b} = \grade{\vec{a} \vec{b}}{0} = \cos{\alpha}$, where $\alpha$ is the angle between the two lines ($\pm1$ when they coincide or are parallel), while $\vec{a}\wedge \vec{b} = \grade{ \vec{a}  \vec{b}}{2}$ is their  intersection point.  If we call the normalized intersection point  $\vec{P}$  (using the appropriate norm), then $\grade{ \vec{a}  \vec{b}}{2} = (\sin{\alpha})\vec{P}$ when the lines intersect and $\grade{ \vec{a}  \vec{b}}{2} = d_{ab}\vec{P}$  when the lines are parallel and are separated by a distance $d_{ab}$. Here we see the remarkable functional polymorphism mentioned earlier, due to the coordination of the two norms.

\item \textbf{Product of two euclidean points.}
\begin{align*}
\vec{P}\vec{Q} = \grade{\vec{P} \vec{Q}}{0} + \grade{\vec{P} \vec{Q}}{2} = -1 + d_{PQ}\vec{V}
\end{align*} 
The inner product (grade-0 part) of any two normalized euclidean points is -1, while the grade-2 part is the direction perpendicular to the joining line $\vec{P} \vee \vec{Q}$.  We sometimes write $\grade{\vec{P} \vec{Q}}{2}$ as $\vec{P} \times \vec{Q}$.  Call $\vec{V}$ the normalized form of $\vec{P} \times \vec{Q}$, then the formula shows that the distance $d_{PQ}$ between the two points satisfies $d_{PQ} = \|\vec{P} \times \vec{Q} \|_\infty$.  We see that the inner product of two points cannot be used to obtain their distance but the grade-2 part can.

\item \textbf{Product of euclidean point and euclidean line.}

This yields a line and a pseudoscalar, both of which contain important geometric information:
\myvspace{-.05in}
\begin{align*}
\vec{a}\vec{P} = \grade{\vec{a} \vec{P}}{1} + \grade{\vec{a} \vec{P}}{3} &= \vec{a} \cdot \vec{P} + \vec{a} \wedge \vec{P} \\
&= \vec{a}_P^\perp + d_{aP}\eye
\end{align*} 
Here $\vec{a}_P^\perp:=\vec{a} \cdot \vec{P}$ is the line passing through $\vec{P}$ perpendicular to $\vec{a}$, 
while the pseudoscalar part has weight $d_{aP}$, the euclidean distance between the point and the line. Note that this inner product is anti-symmetric: $\vec{P} \cdot \vec{a} = - \vec{a} \cdot \vec{P}$.

\end{compactitem}

\myboldhead{Remark} You might be wondering, why is $\vec{a}\cdot\vec{P}$  a line through $\vec{P}$ perpendicular to $\vec{a}$? This is a good opportunity to practice thinking in duality.  Consider $\vec{P}$ as the set of all lines passing through $\vec{P}$ (called the \emph{line pencil} in $\vec{P}$), just as we consider (dually) a line to consist of all the points lying on it.  Indeed,  in the dual exterior algebra where we are operating, $\vec{P}$ -- as a 2-vector --  \textbf{is} composed of 1-vectors (lines) just in this way.  Taking the inner product of $\vec{P}$ with the line $\vec{a}$ removes $\vec{a}$ from the pencil in $\vec{P}$. (That's why the inner product is sometimes called a \emph{contraction} since it reduces the dimension.) The effect is to remove the line  through $\vec{P}$ parallel to $\vec{a}$.  When this line is \emph{completely} removed from $\vec{P}$ it leaves the line  perpendicular to $\vec{a}$. 
 
In the above results, you can also allow one or both of the arguments to be ideal; one obtains in all cases meaningful, \quot{polymorphic} results that in the interests of space we omit. Interested readers can consult \cite{gunn2017b}.  We collect a sample of these formulas in \Tab{tab:pga2d}.  Note that the formulas assume \emph{normalized} arguments. 

After this brief excursion into the world 2-way products, we turn our attention to 3-way products with a repeated factor.  First, we look at products of the form $\vec{XXY}$ (where $\vec{X}$ and $\vec{Y}$ are either $1$- or $2$-vectors). Applying the associativity of the geometric product produces \quot{formula factories}, yielding  a wide variety of important geometric identities. Secondly,  products of the form $\vec{aba}$ for 1-vectors $\vec{a}$ and $\vec{b}$ are used to develop an elegant representation of euclidean motions in PGA based on so-called \emph{sandwich} operators. 
\myvspace{-.2in}
\subsection{Formula factories through associativity} 
  
First recall that for a normalized euclidean point or line, $\vec{X}^2=\pm1$. Use this and associativity to write 
\[ \vec{Y} = \pm(\vec{X}\vec{X})\vec{Y} = \pm\vec{X}(\vec{X}\vec{Y})\] where $\vec{Y}$ is also a normalized euclidean $1$- or $2$-vector.
 The right-hand side yields an \emph{orthogonal decomposition} of $\vec{Y}$ in terms of $\vec{X}$.   Associativity of the geometric product shows itself here to be a powerful tool.  These decompositions are not only useful in their own right, they provide the basis for a family of other constructions, for example,  \quot{the point on a given line closest to a given point}, or \quot{the line through a given point parallel to a given line} (see also \Tab{tab:pga2d}).
 
 Note that the grade of the two vectors can differ.   We work out below three orthogonal projections. 
 As in the above discussions, we  assume the points and lines on the right-hand side of the results represent normalized points and lines, so their coefficients carry unambiguous metric information.  
 
 \begin{figure}[t]
  \centering
\def\xyq{.45}
{\setlength\fboxsep{0pt}\fbox{\includegraphics[width=\xyq\columnwidth]{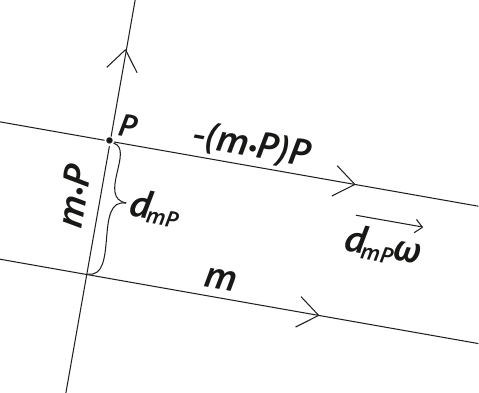}}}\hspace{.005in}
{\setlength\fboxsep{0pt}\fbox{\includegraphics[width=\xyq\columnwidth]{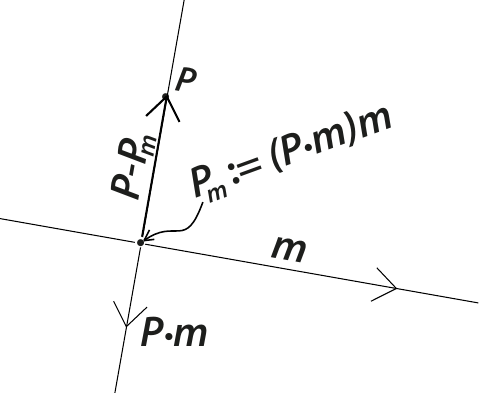}}}\\ \vspace{.007in}
{\setlength\fboxsep{0pt}\fbox{\includegraphics[width=.92\columnwidth]{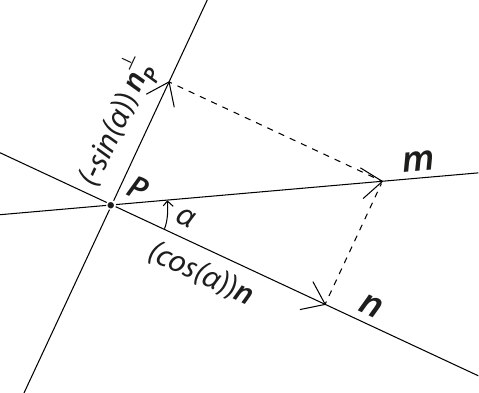}}}
\caption{Orthogonal projections (ul. to lr.): line $\vec{m}$ onto point $\vec{P}$, point  $\vec{P}$  onto line $\vec{m}$, and line $\vec{m}$ onto line $\vec{n}$.}
\label{fig:orthProj}
\end{figure}


\myboldhead{Project line onto line}
Multiply \[\vec{m}\vec{n} = \vec{m}\cdot \vec{n} + \vec{m} \wedge \vec{n}\] with $\vec{n}$ on the right and use $\vec{n}^{2}=1$ to obtain 
\begin{align*}
\vec{m} &= (\vec{m} \cdot \vec{n})\vec{n} + (\vec{m} \wedge \vec{n}) \vec{n} \\
&= (\cos{\alpha}) \vec{n} + (\sin{\alpha} )\vec{P}\vec{n} \\
&= (\cos{\alpha}) \vec{n} - (\sin{\alpha}) \vec{n}^{\perp}_{\vec{P}}
\end{align*}
\myvspace{-.2in}
Thus one obtains a decomposition of $\vec{m}$ as the linear combination of $\vec{n}$ and the perpendicular line $\vec{n}^{\perp}_{\vec{P}}$ through $\vec{P}$.  See \Fig{fig:orthProj}, bottom.


\myboldhead{Project line onto point}
Multiply \[\vec{m}\vec{P} = \vec{m}\cdot \vec{P} + \vec{m} \wedge \vec{P}\] with $\vec{P}$ on the right and use $\vec{P}^{2}=-1$ to obtain 
\begin{align*}
\vec{m} &= -(\vec{m} \cdot \vec{P})\vec{P} - (\vec{m} \wedge \vec{P}) \vec{P} \\
&= -\vec{m}^{\perp}_{\vec{P}}\vec{P} - (d_{\vec{m}\vec{P}}\eye) \vec{P} \\
&= \vec{m}^{||}_{\vec{P}} - d_{\vec{m}\vec{P}} \vec{\omega}
\end{align*}
In the third equation, $\vec{m}^{||}_{\vec{P}}$ is the line through $\vec{P}$ parallel to $\vec{m}$, with the same orientation.  Thus one obtains a decomposition of $\vec{m}$ as the sum of a line through $\vec{P}$ parallel to $\vec{m}$ and a multiple of the ideal line.  Note that just as adding an ideal point (\quot{vector}) to a point translates the point, adding an ideal line to a line translates the line.  
See \Fig{fig:orthProj}, upper left.


\myboldhead{Project point onto line}
Finally one can project a point $\vec{P}$ onto a line $\vec{m}$.  One obtains thereby a decomposition of $\vec{P}$ as $\vec{P}_m$, the point on $\vec{m}$ closest to $\vec{P}$,   plus a vector perpendicular to $\vec{m}$. See \Fig{fig:orthProj}, upper right. Details are in \cite{gunn2017b}.

\myvspace{-.05in}

\subsection{Representing isometries as sandwiches}
\label{sec:isom}
  \begin{figure}[h]
   \centering
   \def\xyz{.95}
{\setlength\fboxsep{0pt}\fbox{\includegraphics[width=\xyz\columnwidth]{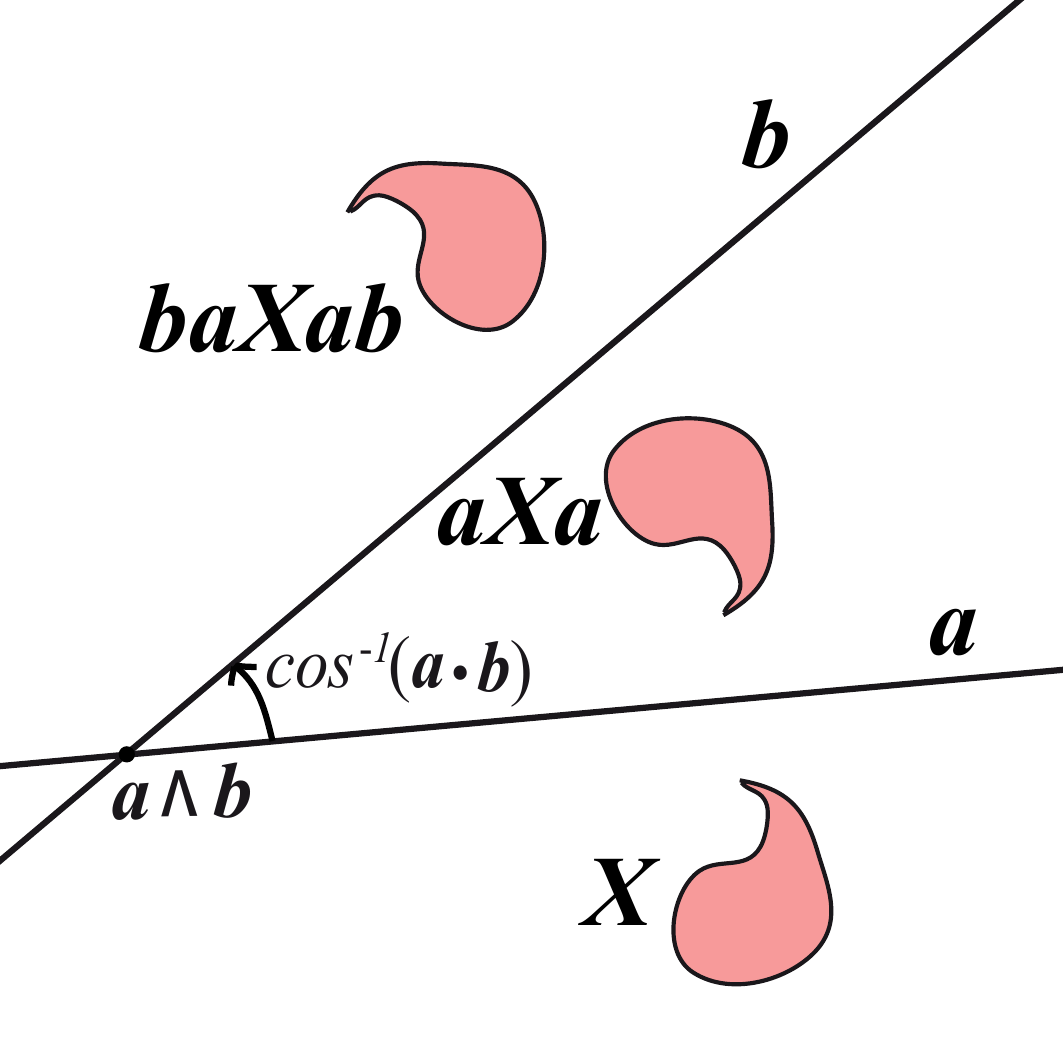}}}
\caption{The reflection in the line $\vec{a}$ is implemented by the sandwich $\vec{a}\vec{X}\vec{a}$; the product of the reflection in line $\vec{a}$ followed by reflection in (non-parallel) line $\vec{b}$ is a rotation around their common point $\vec{a}\wedge\vec{b}$ through $2 \cos^{-1}(\vec{a}\cdot\vec{b})$.}
\label{fig:rotationab}
\end{figure}

Three-way products of the form $\vec{a}\vec{b}\vec{a}$ for euclidean 1-vectors $\vec{a}$ and $\vec{b}$ turn out to represent the reflection of the line $\vec{b}$ in the line $\vec{a}$, and form the basis for an elegant realization of euclidean motions as sandwich operators.  We sketch this here. 

Let $\vec{a}$ and $\vec{b}$ be normalized 1-vectors representing different lines. Then 
\begin{align*}
\vec{a}\vec{b}\vec{a} &= (\vec{a}\vec{b})\vec{a} = (\vec{b} \cdot \vec{a} - \vec{b} \wedge \vec{a}) \vec{a} \\
&= \cos(\alpha)  \vec{a} - \sin(\alpha)\vec{P} \vec{a} \\
&= \cos(\alpha)  \vec{a} +\sin(\alpha) \vec{a}^\perp_\vec{P}
\end{align*}
Compare that with the orthogonal decomposition for $\vec{b}$ obtained above: $\vec{b}=\cos(\alpha) \vec{a} - \sin(\alpha) \vec{a}^\perp_\vec{P}$.
Using the fact that $\vec{a}^\perp_\vec{P}$ is a line perpendicular to $\vec{a}$ leads to the conclusion that $\vec{a}\vec{b}\vec{a}$ must be the reflection of $\vec{b}$ in $\vec{a}$, since the reflection in $\vec{a}$ is the unique linear map fixing $\vec{a}$ and mapping $\vec{a}^\perp_\vec{P}$ to $-\vec{a}^\perp_\vec{P}$. We call this the \emph{sandwich operator} corresponding to $\vec{a}$ since $\vec{a}$ appears on both sides of the expression. It's not hard to show that for a euclidean point $\vec{P}$, $\vec{a}\vec{P}\vec{a}$ is the reflection of $\vec{P}$ in the line $\vec{a}$ . Similar results apply in higher dimensions: the same sandwich form for a reflection works regardless of the grade of the \quot{meat} of the sandwich.

\myboldhead{Rotations and translations} It is well-known that all isometries of euclidean space are generated by reflections.   The sandwich $\vec{b}(\vec{a}\vec{X}\vec{a})\vec{b}$ represents first reflecting in line $\vec{a}$, then reflecting in line $\vec{b}$.  When the lines meet at angle $\frac{\alpha}2$, this is well-known to be  a rotation around the point $\vec{P}$ through of angle $\alpha$. Writing $\vec{R} := \vec{a}\vec{b}$, the  rotation can be expressed as $\vec{R} \vec{X} \widetilde{\vec{R}}$. (Here, $\widetilde{\vec{R}}$ is the \emph{reversal} of $\vec{R}$, obtained by writing all products in the reverse order). 
See \Fig{fig:rotationab}.  When $\vec{a}$ and $\vec{b}$ are parallel, $\vec{R}$ generates the \emph{translation} in the direction perpendicular to the two lines, of twice the distance between them. 
When $\vec{R}$ is normalized so that $\vec{R}\widetilde{\vec{R}} = 1$, it's called a \emph{rotor}.  

\myboldhead{Exponential form for rotors} Rotors can be generated directly from the normalized center point $\vec{P}$ and angle of rotation $\alpha$ using the exponential form $\vec{R} = e^{\frac{\alpha}{2}\vec{P}}$.  This is another standard technique in geometric algebra: The exponential behaves like the exponential of a complex number since a normalized euclidean point satisfies $\vec{P}^2=-1$.  When $\vec{P}$ is ideal ($\vec{P}^2=0$), the same process yields a translation through distance $d$ perpendicular to the direction of $\vec{P}$, by means of the formula $\vec{T} = e^{\frac{d}{2}\vec{P}} = 1 + \frac{d}{2}\vec{P}$. 

\Tab{tab:pga2d} contains an overview of formulas available in $\pdclal{2}{0}{1}$, most of which have been introduced in the above discussions. We are not aware of any other frameworks offering comparably concise and polymorphic formulas for plane geometry. 


\begin{table}
\begin{center}
\scriptsize 
\renewcommand{\arraystretch}{1.15}
\begin{tabular}[h]{||  l  |    c ||}  \hline
Operation  & PGA \\ \hline  \hline
Intersection point of two lines &   $ \vec{a} \wedge \vec{b}$  \\  \hline
Angle of two intersecting lines &  $ \cos^{-1}(\vec{a} \cdot \vec{b} ) $ \\ \cline{2-2}
&  $\sin^{-1}(\|\vec{a}\wedge\vec{b}\|) $ \\ \hline
Distance of two $||$ lines &  $\| \vec{a} \wedge \vec{b}\|_\infty $ \\ \hline
Joining line of two points  &     $\vec{P} \vee \vec{Q}$ \\ \hline
$\perp$ direction to join of two points  &     $\vec{P} \times \vec{Q}$ \\ \hline
Distance between two  points &  $\|\vec{P}\vee \vec{Q}\|$,$~~\| \vec{P} \times \vec{Q} \|_\infty $ \\ \hline
Oriented distance  point to line  &  $\| \vec{a} \wedge \vec{P} \| $  \\ \hline
Angle of ideal point to line  &  $\sin^{-1}{(\| \vec{a} \wedge \vec{P} \|_\infty)} $  \\ \hline
Line through point $\perp$ to line &   $\vec{P} \cdot \vec{a} $ \\ \hline
Nearest point on line to point &   $(\vec{P} \cdot \vec{a}) \vec{a}$ \\ \hline
Line through point $||$ to line &   $(\vec{P} \cdot \vec{a})  \vec{P}$ \\ \hline
Area of triangle $ABC$ & $\frac12 \|\vec{A} \vee \vec{B} \vee \vec{C} \| $ \\ \hline
Reflection in line ($\vec{X} =$ point or line) & $\vec{a} \vec{X} \vec{a}$ \\ \hline \renewcommand{\arraystretch}{1.3} 
Rotation around point of angle $2 \alpha$ & $\vec{R} \vec{X} \widetilde{\vec{R}} ~~ (\vec{R} :=  e^{\alpha \vec{P}}$) \\ \hline
Translation by $2d$ in direction $\vec{V}^\perp$& $\vec{T} \vec{X} \widetilde{\vec{T}} ~~ (\vec{T} := 1 + d{\vec{V}})$  \\ \hline
\end{tabular}
\end{center}
\myvspace{.1in}
\caption{A sample of geometric constructions and formulas in the euclidean plane using PGA (assuming normalized arguments, all arguments euclidean unless otherwise stated). }
\label{tab:pga2d}
\end{table}

\subsection{Automatic differentiation}
\label{sec:ad}
\cite{hessob87} introduces the term \quot{geometric calculus} for the application of calculus to geometric algebras, and shows that it offers an attractive unifying framework in which many diverse results of calculus and differential geometry can be integrated.  While a treatment of geometric calculus lies outside the scope of this article (see \cite{dfm07}, Ch. 8, for a practical introduction), we want to present a related result to give a flavor of what is possible in this direction. 

Notice that the elements $\{1,\eye\}$ generate a 2-dimensional sub-algebra of $\pdclal{n}{0}{1}$ consisting of scalars and pseudoscalars.  This algebra is known as the \emph{dual numbers} and can be abstractly characterized by the fact that $1^2=1$ while $\eye^2=0$.  Already Eduard Study, the inventor of dual numbers, realized that they can be used to do automatic differentiation  (\cite{study03}, Part II, \S 23).
A modern reference describes how \cite{wikiAD}:
\begin{quote}
Forward mode automatic differentiation is accomplished by augmenting the algebra of real numbers and obtaining a new arithmetic. An additional component is added to every number which will represent the derivative of a function at the number, and all arithmetic operators are extended for the augmented algebra. The augmented algebra is the algebra of dual numbers.
\end{quote}
This extension can be obtained by beginning with the monomials. Given $p_k(x)=x^k$, define \[p_k(x+y\eye) := (x+y\eye)^k = x^k+nx^{n-1}y\eye\] All higher terms disappear since $\eye^2=0$.  Setting $y=1$ we obtain \[{p}_k(x+\eye) = p_k(x)+\dot{p}_k(x)\eye\] That is, the scalar part is the original polynomial and the pseudoscalar, or dual, part is its derivative. In general if $u$ is a function $u(x)$ with derivative $\dot{u}$, then \[p_k(u+\dot{u}\eye) = p_k(u) + \dot{p}_k(u)\eye\] Thus, the coefficient of $\eye$ tracks the derivative of $p_k$.  Extend these definitions to polynomials by additivity in the obvious way.  Since the polynomials are dense in the analytic functions, the same \quot{dualization} can be extended to them and one obtains in this way robust, exact automatic differentiation.  One can also handle multivariable functions of $n$ variables, using the $(n)$ ideal $n-$vectors $E_i$ for $i>0$ (representing the ideal directions of euclidean $n$-space) as the nilpotent elements instead of $\eye$. 
For a live JavaScript demo see \cite{ganjacs}.

\myvspace{-.2in}
\section{PGA for euclidean space: $\pdclal{3}{0}{1}$}
\label{sec:eucspace}
If you have followed the treatment of plane geometry using PGA, then you are well-prepared to tackle the 3D version $\pdclal{3}{0}{1}$.  
Naturally in 3D one has points, lines, \textbf{and} planes, with the planes taking over the role of lines in 2D (as dual to points); the lines represent a new, middle element not present in 2D. 
A look at the table of formulas for 3D (\Tab{tab:pga3d}) confirms that many of the 2D formulas reappear, with planes substituting for lines.
If you re-read Examples \ref{sec:3dkal} and  \ref{sec:3dscrew} now you should understand much better how 3D isometries are represented in PGA, based on what you've learned about 2D sandwiches.

In the interests of space, we leave it to the reader to confirm the similarities of the 3D case to the 2D case and instead focus for the remainder of this article on one important difference:   bivectors of $\pdclal{3}{0}{1}$, which, as we mentioned above, have no direct analogy in $\pdclal{2}{0}{1}$.  

\begin{table}[t]
\begin{center}
\scriptsize 
\renewcommand{\arraystretch}{1.15}
\begin{tabular}[h]{||  l  |    c ||}  \hline
\textbf{Operation}  & \textbf{formula} \\ \hline  \hline
Intersection line of two planes &   $ \vec{a} \wedge \vec{b}$  \\  \hline
Angle of two intersecting planes &  $ \cos^{-1}(\vec{a} \cdot \vec{b} ) $ \\ \cline{2-2}
 &  $\sin^{-1}(\|\vec{a}\wedge\vec{b}\|) $ \\ \hline
Distance of two $||$ planes &  $\| \vec{a} \wedge \vec{b} \|_\infty$ \\ \hline
Joining line of two points  &     $\vec{P} \vee \vec{Q}$\\ \hline
Intersection point of three planes &    $\vec{a} \wedge \vec{b} \wedge \vec{c}$ \\ \hline
Joining plane of three points &    $\vec{P} \vee \vec{Q} \vee \vec{R}$ \\ \hline
Distance from  point to plane  &  $\| \vec{a} \wedge \vec{P} \|$  \\ \hline
Angle of ideal point to plane  &  $\sin^{-1}{(\| \vec{a} \wedge \vec{P} \|_\infty)}$  \\ \hline
$\perp$ line to join  of two  points &  $ \vec{P} \times \vec{Q}$ \\ \hline
Distance of two points &  $\|\vec{P}\vee \vec{Q}\|$,$~~\| \vec{P} \times \vec{Q} \|_\infty $\\ \hline
Line through point $\perp$ to plane  &   $\vec{P} \cdot \vec{a} $ \\ \hline
Closest point on plane to point &   $(\vec{P} \cdot \vec{a}) \vec{a}$ \\ \hline
Plane through point $||$ to plane &   $(\vec{P} \cdot \vec{a})  \vec{P}$ \\ \hline
Plane through line $\perp$ to plane  &   $\velo \cdot \vec{a} $ \\ \hline
Intersection of line and plane & $\velo \wedge \vec{a}$ \\ \hline
Joining plane of point and line &   $\vec{P} \vee \velo $ \\ \hline
Plane through point $\perp$ to line  & $ \vec{P} \cdot  \velo $ \\ \hline
Closest point on line to point & $ ( \vec{P} \cdot \velo )  \velo $ \\ \hline
Line through point $||$ to line &   $(\vec{P} \cdot \velo)  \vec{P}$ \\ \hline
Line through point $\perp$ to line &   $((\vec{P} \cdot \velo)  \velo)\vee \vec{P}$ \\ \hline
Volume of tetrahedron $ABCD$ & $\frac16 \|\vec{A} \vee \vec{B} \vee \vec{C}\vee\vec{D} \| $ \\ \hline
Refl. in plane ($\vec{X} =$ pt, ln, or pl)  & $\vec{a} \vec{X} \vec{a}$ \\ \hline \renewcommand{\arraystretch}{1.4}
Rotation with axis $\velo$ by angle $2 \alpha$ & $\vec{R} \vec{X} \widetilde{\vec{R}} ~~ (\vec{R} :=  e^{\alpha \velo}$) \\ \hline
Translation by $2d$ in direction $\vec{V}$ & $\vec{T} \vec{X} \widetilde{\vec{T}} ~~ (\vec{T} := {(\EE{0} \vee \vec{V})\eye})$  \\ \hline
Screw with axis $\velo$ and pitch $p$ &$ \vec{S} \vec{X} \widetilde{\vec{S}} ~~ (\vec{S} := e^{t(1+p\eye)\velo})$  \\ \hline
\end{tabular}
\end{center}
\myvspace{.1in}
\caption{A sample of geometric constructions and formulas in 3D using PGA (assuming normalized arguments).}
\label{tab:pga3d}
\end{table}

\myvspace{-.2in}

\subsection{Lines and 2-vectors in 3D}
\label{sec:l2v2d}
In $\pdclal{2}{0}{1}$, all $k$-vectors are \emph{simple}, that is, they can be written as the product of $k$ 1-vectors.  This is no longer the case in $\pdclal{3}{0}{1}$.  
If $\sigo_1$ and $\sigo_2$ are bivectors representing {skew} lines then the sum $\momo := \sigo_1 + \sigo_2$ is non-simple.  
(For a proof, see \cite{gunn2011} or \cite{gunnthesis}, Ch. 8.) Non-simple bivectors make the 3D case much complex and interesting than the 2D case.  Due to space limitations we can only give a flavor of this behavior in what follows.

\begin{figure}[t]
   \centering\
   \def\xyy{.95}
   \def\xyz{.5} 
   \def\xyw{.39} 
{\setlength\fboxsep{0pt}\fbox{\includegraphics[width=\xyz\columnwidth,trim={5mm 0 9mm 0}, clip]{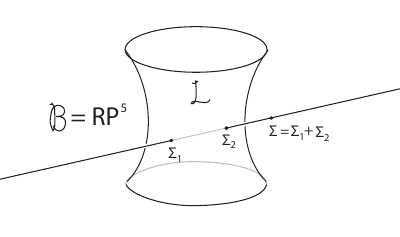}}} \hspace{.04in}  
{\setlength\fboxsep{0pt}\fbox{\includegraphics[width=\xyw\columnwidth,trim={6cm 0 2cm 0},clip]{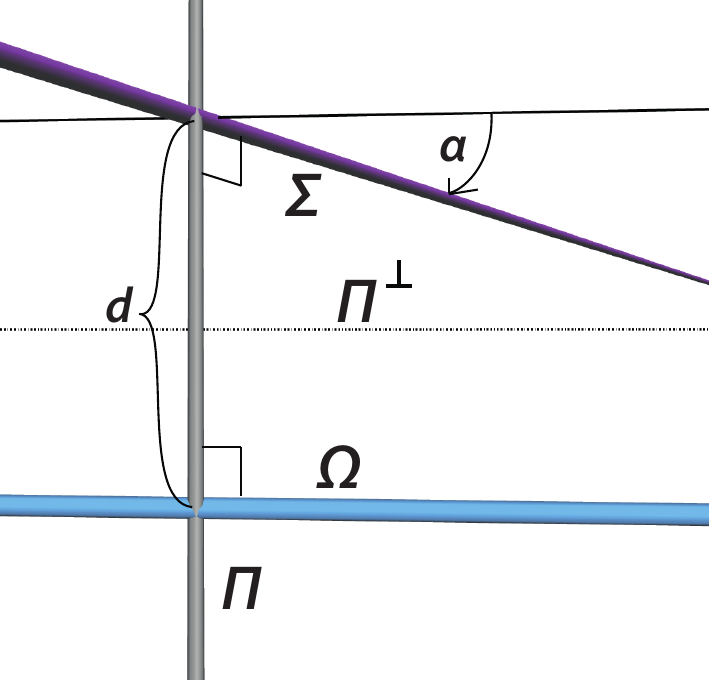}}}
\caption{\emph{Left:}The space of lines sits inside the space of 2-vectors as a quadric surface $\mathcal{L}$. \emph{Right}: Product of two skew lines $\velo$ and $\sigo$ involving the common normals $\momo$ (euclidean) and $\momo^\perp$ (ideal).}
\label{fig:complexspace}
\end{figure}

 The space of bivectors is  spanned by the 6 basis elements $\e{ij} := \e{i} \e{j}$.  The $\e{ij}$ can be thought of as the lines of intersection of the 4 basis planes. The three elements $\e{0i}$ are ideal lines 
 while $(\e{23}, \e{31}, \e{12})$ are lines through the origin in the $(x, y, z)$-directions, resp. The bivectors form a 5-dimensional projective space $\mathcal{B} := \RP{5}$. The condition that a bivector is simple (represents a line) can be translated into a quadratic constraint on the coordinates of the bivector a result due to  Pl\"{u}cker (1801-1868).  This defines the \emph{Pl\"{u}cker quadric} $\mathcal{L}$, a 4D surface sitting inside $\mathcal{B}$, associated to the well-known \emph{Pl\"{u}cker coordinates} for lines. Points not on the quadric are non-simple bivectors, also known as \emph{linear line complexes}. Consult Figure \ref{fig:complexspace}. 

\myvspace{-.15in}
\subsection{Product of two euclidean lines}
\label{sec:prtwoli}
Due to space limitations we can only give a small taste of 3D line geometry, by calculating the geometric product of two lines.
Let the two lines be $\velo$ and $\sigo$.  Assume they are euclidean, skew, and normalized, i.e., $\velo \wedge \sigo \neq 0$ and $\velo^2 = \sigo^2 = -1$. Two euclidean lines determine in general a unique third euclidean line that is perpendicular to both, call it $\momo$.   Consult \Fig{fig:complexspace}, right. 
$\velo\sigo$ consists of 3 parts, of grades 0, 2, and 4:
\begin{align*}
\velo\sigo &= \grade{\velo\sigo}{0} + \grade{\velo\sigo}{2} + \grade{\velo\sigo}{4} \\
&= \velo \cdot \sigo + \velo \times \sigo + \velo \wedge \sigo \\
&= \cos{\alpha} + (\sin{\alpha}\momo + d \cos{\alpha}\momo^\perp) + d \sin{\alpha} \eye 
\end{align*}
 Here $\alpha$ is the angle between $\velo$ and $\sigo$, viewed along the common normal $\momo$; $d$ is the distance between the two lines measured along $\momo$. $d \sin{\alpha}$ is the volume of a tetrahedron determined by unit length segments on $\velo$ and $\sigo$. Finally, $ \velo \times \sigo$ is a weighted sum of  $\momo$ and $\momo^\perp$. The appearance of $\momo^\perp$ is not so surprising,  as it is also a \quot{common normal} to $\velo$ and $\sigo$, but as an ideal line, is easily overlooked. 
 

Does $\velo \sigo$ have a geometric meaning? 
Consider sandwich operators with bivectors, that is, products of the form $\velo\vec{X}\widetilde{\velo} $ 
for simple $\velo$. Such a product is called a \emph{turn} since it implements a half-turn around the axis $\velo$.  And, in turn, the turns generate the full group $\Eucgd{3}$ of direct euclidean isometries (\cite{study91}. A little reflection shows that the  composition of the two turns $\velo\sigo$ will be a screw motion that rotates around the common normal $\momo$ by $2\alpha$ while translating by $2d$ in the same direction (the latter is a \quot{rotation} around $\momo^\perp$).   This is analogous to the product of two reflections meeting at angle $\alpha$ discussed above in \Sec{sec:isom}. 

\subsection{Remarks on Kinematics and Mechanics}
Hopefully the preceding remarks make clear the central role that bivectors play in 3D PGA. Essentially they form the \emph{Lie algebra} of $\Eucgd{3}$, the oriented Lie group of euclidean space $\Euc{3}$.  Exponentiating them produces the \emph{Spin} group  in the even subalgebra $\pdclplus{3}{0}{1}$ consisting of the normalized  elements ($\vec{G}\widetilde{\vec{G}}=1$).  These elements are sometimes called \emph{rotors}, and form  a 2:1 covering of  $\Eucgd{3}$.  This even sub-algebra in turn is isomorphic to the dual quaternions.

When one applies this framework to calculating  the motion of the free top one obtains the following Euler equations of motion:
\begin{align}\label{eqn:eqnmot}
\dot{\vec{g}} &= \vec{g}\velo_c\\
\dot{\velo}_c &= 2\inert^{-1}(\inert(\velo_c) \times \velo_c )
\end{align}
Here $\vec{g}(t)$ is a path in the  Lie group representing the motion of the body.  $\vec{A}$ is the inertia tensor  of the body.  $\velo_c$ (\emph{resp.}, $\inert(\velo_c)$) is the velocity of the body (\emph{resp.}, momentum of the body) in body coordinates; both are represented by bivectors in $\pdclal{3}{0}{1}$. See \cite{gunn2011} or \cite{gunnthesis}, Ch. 9, for details.  For a very compact PGA implementation see \cite{ganjacs}.  

These equations behave particularly well numerically: the solution space has 12 dimensions (the isometry group is 6D and the momentum space (bivectors) also) while the integration space has 14 dimensions ($\pdclplus{3}{0}{1}$ has dimension 8 and the space of bivectors has 6).  Normalizing the computed rotor $\vec{g}$ brings one directly back to the solution space.  In traditional matrix approaches as well as the CGA approach  (\cite{lasenby2011}), the co-dimension of the solution space within the integration space is much higher and leads typically to the use Lagrange multipliers or similar methods to maintain accuracy.  This advantage over VLAAG and CGA is typical of the PGA approach for many related computing challenges.

 \myvspace{-.15in}
 \section{Implementation}
 Our description would be incomplete without discussion of the practical issues of implementation.  This has  been the focus of much work and there exists a well-developed theory and practice for general geometric algebra implementations to maintain performance parity with traditional approaches.  See \cite{hildebrand13}.  PGA presents no  special challenges in this regard; in fact, it demonstrates clear advantages over other geometric algebra approaches to euclidean geometry in this regard (\cite{gunn2017a}). 
For a full implementation of PGA in JavaScript ES6 see Steven De Keninck's ganja.js project on GitHub \cite{ganja} and the interactive example set at \cite{ganjacs}.
 \begin{table}
\begin{center}
\scriptsize
\setlength{\extrarowheight}{4pt}
\def\xyz{3pt}
\begin{tabular}{| p{.45\columnwidth} | p{.45\columnwidth} |}
\hline
\textbf{PGA} & \textbf{VLAAG} \\
\hline
\hline
Unified representation for points, lines, and planes based on a graded exterior algebra; all are \quot{equal citizens} in the algebra.
&
The basic primitives are  points and vectors and all other primitives are built up from these. For example, lines in 3D sometimes  parametric, sometimes w/ Pl\"{u}cker coordinates. \\[\xyz]
\hline
Projective exterior algebra provides robust meet and join operators that deal correctly with parallel entities.
&
Meet and join operators only possible when homogeneous coordinates are used, even then tend to be \emph{ad hoc} since points have distinguished role and ideal elements  rarely integrated. \\[\xyz]
\hline
Unified, high-level treatment of euclidean (\quot{finite}) and ideal (\quot{infinite}) elements of all dimensions. Unifies e.g. rotations and translations, simple forces and force couples.
&
Points (euclidean) and vectors (ideal) have their own rules, user must keep track of which is which; no higher-dimensional analogues for lines and planes.\\[\xyz]
\hline
Unified repn. of isometries based on sandwich operators which act uniformly on points, lines, and planes.
&
Matrix representation for isometries has different forms for points, lines, and planes.\\[\xyz]
\hline
Same representation for operator and operand: $\vec{m}$ is the  plane as well as the reflection in the plane.
&
Matrix representation for reflection in $\vec{m}$ is different from the vector representing the plane. \\[\xyz]
\hline
Compact, universal expressive formulas and constructions based on geometric product (see Tables \ref{tab:pga2d} and \ref{tab:pga3d}) valid for wide range of argument types and dimensions.
&
Formulas and constructions are \emph{ad hoc}, complicated, many special cases, separate formulas for points/lines/planes, for example, compare   \cite{gg90}. \\[\xyz]
\hline
Well-developed theory of implementation optimizations to maintain performance parity.
&
Highly-optimized libraries, direct mapping to current GPU design. \\[\xyz]
\hline
Automatic differentiation of real-valued functions using dual numbers. 
&
Numerical differentiation \\
\hline
\end{tabular}
\myvspace{.1in}
\caption{A comparison of PGA with the standard VLAAG approach. }
\label{tab:comp}
\end{center}
\end{table}

\section{Comparison}
\label{sec:comp}
Table \ref{tab:comp} encapsulates the foregoing results in a feature-by-feature comparison with the standard (VLAAG) approach.   It establishes that PGA fulfills all the features on our wish-list in Sec. \ref{sec:wishlist}, while the standard approach offers almost none of them.  (For a proof that PGA is coordinate-free, see the Appendix in \cite{gunn2017b}.) 

\subsection{Conceptual differences}
How can we characterize conceptually the difference of the two approaches leading to such divergent results? First and foremost:  VLAAG is \emph{point-centric}: other geometric primitives of VLAAG such as lines and planes are built up out of points and vectors. PGA on the other hand is \emph{primitive-neutral}: the exterior algebra(s) at its base provide \emph{native} support for the subspace lattice of points, lines and planes (with respect to both join and meet operators). Secondly, the projective basis of PGA allows it to deal with points and vectors in a unified way: vectors are just ideal points, and in general, the ideal elements play a crucial role in PGA to integrate parallelism, which typically has to be treated separately in VLAAG.  The existence of the ideal norm in PGA goes beyond the purely projective treatment of incidence, producing polymorphic metric formulas that, for example, correctly handle two intersecting lines whether they intersect or are parallel (see above \Sec{sec:prodpr}).  We believe that the resulting tables of formulas (\Tab{tab:pga2d} and \Tab{tab:pga3d}), based on this \quot{dynamic duo} of standard and ideal norms,  are \quot{world champions} with respect to compactness and polymorphicity among all existing frameworks for euclidean geometry, and that there are many more formulas waiting to be discovered (after all, we've only considered the 2-way products and a small subset of the 3-way products). Compare  \cite{gg90} for selected VLAAG analogs. The representation of isometries using sandwich operators generated by reflections in planes (or lines in 2D) can be understood as a special case of this \quot{compact polymorphicity}: the sandwich operator $\vec{g}\vec{X}\tilde{\vec{g}}$ works no matter what $X$ is, the same representation works whether it appears as operator or as operand, and rotations and translations are handled in the same way.

\subsection{The universality of PGA}
The previous section has looked for and found the basis for the superiority of PGA over VLAAG in its structural basis. Here we go further and show that a large extent, alternate approaches to euclidean geometry are present already in PGA as parts of the whole.

\myboldhead{Vector algebra} The previous section has already suggested that VLAAG can be seen less as a direct competitor to PGA than as a restricted subset.  Indeed, restricting attention to the vector space of $n$-vectors (sometimes written $\bigwedge^n$) in PGA  essentially yields standard vector algebra. Define the \quot{points} to be euclidean $n$-vectors ($\vec{P}^2 \neq 0$) and \quot{vectors} to be ideal $n$-vectors ($\vec{P}^2 = 0$).  All the rules of vector algebra can be then derived using the   vector space structure of $\bigwedge^n$ along with the standard and ideal PGA norms (assuming normalized arguments as usual).  
This embedding of vector algebra in PGA also comes with a nice geometric intuition absent in traditional vector algebra: the vectors make up the ideal plane bounding the euclidean space of points, \emph{ i. e.},  points and vectors make up a connected, unified space. Furthermore, intuitions developed in vector algebra such as \quot{Adding a vector to a point translates the point.} have natural extensions in PGA, since adding an ideal line (plane) to a euclidean line (plane) translates the line (plane) (whereby the two lines must be co-planar).  Such patterns are legion.

\myboldhead{Linear algebra and analytic geometry} Note that PGA is fully compatible with the use of linear algebra -- the difference is that it no longer is needed to implement euclidean motions, a role for which it is not particularly well-suited. In a similar way, we envision the development of an analytic geometry based on the full extent of PGA, not just on the small subset present in VLAAG, and would have at its disposal the geometric calculus sketched in \Sec{sec:ad}.  Traditional analytic geometry would make up a small subset of this extended analytic geometry, like vector algebra within PGA proper.

\myboldhead{Quaternions and dual quaternions} Many  aspects of PGA are present in embryonic form in quaternions and dual quaternions, but they only find their full expression and utility when embedded in the full algebra PGA. Indeed, the quaternion and dual quaternion algebras are isomorphically embedded in the even sub-algebra $\pdclplus{n}{0}{1}$ for $n \ge 3$. 
The advantage of the embedding in PGA are considerable. The full algebraic structure of PGA provides a much richer environment than these quaternion algebras alone.  Few of the formulas in Tables \ref{tab:pga2d} and \ref{tab:pga3d} are available in the quaternion algebras alone since the quaternion algebras only have natural representations for primitives of even grade. For example, in PGA, you can apply the  sandwiches to geometric primitives of any grade.  In contrast, one of the \quot{mysteries} of contemporary dual quaternion usage is that there are separate \emph{ad hoc} representations for points and planes and slightly different forms of the sandwich operator for each in order to be able to apply euclidean isometries.  These eccentricities disappear when, as in PGA, there are native representations for points and planes, see \cite{gunn2017a}, \S 3.8.1. 
The PGA embedding clears up other otherwise mysterious aspects of current dual quaternion practice. Consider the  dual unit $\epsilon$ satisfying $\epsilon^2=0$. In the embedding map, it maps to the pseudoscalar $\eye$ of the algebra (for details see \cite{gunnthesis}, \S 7.6), perhaps tarnishing the mystique  but replacing it with a deeper understanding of the dual quaternions.  


 \section{Migrating to PGA}
The foregoing exposition establishes that PGA is by any metric a strong candidate for doing euclidean geometry on the computer. The natural next question for developers is, what is involved in migrating to PGA from one of the alternatives discussed above?
In fact, the use of homogeneous coordinates and the inclusion of quaternions, dual quaternions, and exterior algebra in PGA means that many practitioners already familiar with these tools can expect a  gentle learning curve. Furthermore, the availability of a JavaScript implementation on GitHub (\cite{ganja}) and the existence of platforms such as Observable notebooks \cite{observable} means that interested users can quickly begin to work and prototype their applications. Readers who would like first to deepen their understanding of the underlying mathematics are referred to the bibliography, particularly \cite{gunn2011} and, for the full metric-neutral treatment, \cite{gunnthesis}.
 
\section{Conclusion} We close with some reflections on the intimate relationship between mathematics and its applications. Naturally there are good reasons to focus on the primacy of the application, and the use of mathematics as a tool to achieve that end. And, indeed, users who persevere in mastering PGA can expect to reap the benefits established in the foregoing discussion.  Hence, existing applications in all the practical fields mentioned at the beginning of this article will, we believe, benefit in this way from exposure to PGA.   Equally exciting in our view is the prospect that PGA, as a new  way of thinking about euclidean geometry, may lead to innovative applications that were, so to speak, hidden from view using previous approaches.  So from whatever direction you are coming to the subject -- whether you are interested in improving an existing application or in learning a new approach to euclidean geometry -- PGA has plenty to offer to all.
 


\ifthenelse{\equal{\isLong}{false}}
{}
{
\section{Further reading}
Readers interested in knowing more about PGA can consult the following publications:
\small
\begin{compactitem}
\item \cite{gunn2011}: initial publication of euclidean PGA, including a self-contained treatment of  kinematics and rigid body mechanics,
\item  \cite{gunnFull2010}: an extended version of \cite{gunn2011},
\item  \cite{gunnthesis}: dissertation featuring a \MN setting for PGA including hyperbolic and elliptic space,
\item \cite{gunn2017a}: comparison of PGA  to CGA (another version of geometric algebra for euclidean geometry).
\item \cite{gunn2017b}: tutorial-like introduction to PGA applied to euclidean plane geometry,
\end{compactitem}
\normalsize
 \section{Application areas}
 We want now to turn to consider its potential impact on three application areas: computer graphics, game engines, and 3D manifolds (a branch of topology) before wrapping up with a comparison with the standard approach to euclidean geometry.
Progress takes various forms.  It can be theoretical or practical.  But the dividing line is often not so clear as it appears.  

 \subsection{Computer graphics}
 
 Here we take \emph{computer graphics} in its wider sense, to include the ecosystem of theory and practice that includes animation, modeling, and rendering.  If one compares the formulas in Table \ref{tab:pga3d} with \cite{gg90}, a similar but smaller collection based on the standard approach, the compactness and elegance of the PGA formulas is striking.\footnote{We would in fact like to issue a challenge to the computer graphics  community: is there another framework which provides these formulas so compactly (measured by symbol count) as PGA,  with the same degree of polymorphicity as PGA?  Candidates should be sent via email to the author.} Consequently, Euclidean PGA shows great promise to provide a unified geometric language for all fields related to computer graphics. Those who see a similar role for CGA are referred to \cite{gunn2017a} discusses this issue in detail.
Having a high-level polymorphic description of euclidean geometry has obvious advantages for problem-solving, since it shortens solutions and reduces the chances for errors.  Computer graphics, however, must also satisfy the constraints of real-time computation.  How does PGA behave with respect to practical implementation?   
  
 \section{Implementation issues}
 Fortunately there has been lots of work carried out to establish that, from a practical perspective, geometric algebras remain competitive with standard linear algebra approaches to the same problems.  See for example Ch. 22 of \cite{dfm07} where a comprehensive strategy is described for avoiding pitfalls in implementing geometric algebra structures. This includes keeping track of the non-zero grades of a multi-vector, automatically generating special purpose code for sub-products of the geometric product involving pure $k$-vectors, and judicious use of linear algebra methods behind the scenes in computationally-intensive loops.  Until GPU technology is based on GA models and not on VLAAG ones, GA toolkits will need to convert from their GA representation to VLAAG API's of the GPU (such as OpenGL or DirectX).  This will of course be hidden from the PGA user/developer.

 \subsection{3D Kinematics via 2-vectors}
\label{sec:kin}
 
 M\"{o}bius, in his investigations of euclidean statics \cite{moebius37}, discovered the linear line complex (which is an older name for the 2-vectors discussed here, naturally in a much less clear form) and thereby initiated the modern study of \emph{line geometry}.  Other investigators, notably Pl\"{u}cker and Felix Klein (\cite{klein71b}), pursued this direction of research and obtained an elegant and complete theory of kinematics and rigid body mechanics based on line geometry (\cite{ziegler}).  In fact, the history of geometric algebra itself is inextricably intertwined with this stream of research:  Clifford was a friend and admirer of Klein's, and his investigations of rigid body mechanics in elliptic space \cite{clifford74} using line geometry cannot be separated from his invention of dual quaternions (\cite{clifford73}) and geometric algebra (\cite{clifford78}).  Consult \cite{ziegler} for a detailed account of these fascinating chapter of $19^{th}$ century mathematics.
 
\subsubsection{The null system}

Since it plays a central role in the following discussion of kinematics and rigid body mechanics, we present M\"{o}bius's key results in \cite{moebius37} in more detail.  He was investigating 

 The basic object of kinematics is to describe the motion of bodies preserving euclidean (or some other metric) measurement.  Such a movement is called an \emph{isometry}.  We describe such a motion as a path $g(t) \subset \Eucgd{3}$, the group of orientation-preserving Euclidean isometries, 
 such that $g(0) = \mathbb{1}$, the identity.  We saw in the example above \Sec{sec:3dscrew} that such a path can be expressed in $\pdclal{3}{0}{1}$ as the rotor $\vec{g}(t) := e^{t(\velo+p\velo^\perp)}$, representing a screw motion with axis $\velo$ and pitch $p$. And, in light of what was said regarding non-simple 2-vectors above, we see that $\velo+p\velo^\perp$ is just a general 2-vector  (with the exception of ideal 2-vectors represented as $p \velo^\perp$). So, the exponential of any 2-vector generates a euclidean motion, and any euclidean motion can be so generated.  Thus, the 2-vectors function as the \emph{Lie algebra} for the group $\Eucg{3}$, and that group itself is neatly imbedded in the even subalgebra of $\pdclal{3}{0}{1}$ (i. e., the subalgebra consisting of even-grade vectors: scalars, bivectors, and pseudoscalars).\footnote{The exponential map actually yields a 2:1 covering map of the direct Euclidean group since the rotors $\vec{R}$ and $-\vec{R}$ yield the same isometry.}
 
 Also the infinitesimal aspects of kinematics are elegantly represented here. For simplicity write $\vec{g}(t) := e^{t \velo}$ for a general bivector $\velo$. Then $\vec{g}(t)$ is a euclidean motion and $\dot{\vec{g}}(t) = \velo \vec{g}(t)$, so $\dot{\vec{g}}(0) = \velo$. That is, the derivative of the motion at time $t=0$ is $\velo$ itself. We call $\velo$ in this context the \emph{global velocity state}. Also forces and momenta are represented by 2-vectors.  For example, a simple force is represented by the line $\momo$ upon which it acts.   
The following theorem shows how $\velo$ provides the basis for calculating all kinematic quantities in PGA. Consult \cite{gunnthesis}, \S 8.3 for a proof.  The subscripts $c$ and $s$ represent the body and space coordinate systems (space limitations prohibit a fuller discussion of this essential element of kinematics). 

\begin{theorem} \label{thm:liebracket}
For  time-varying multi-vector $\vec{X}$ subject to the motion $\vec{g}$ with velocity in the body $\velo_c$,
\[
\dot{\vec{X}}_s = \vec{g}(\dot{\vec{X}}_c+2(\velo_c \times \vec{X}_c )) \tilde{\vec{g}} = \vec{g}\dot{\vec{X}}_c\tilde{\vec{g}} +2(\velo_s\times \vec{X}_s)
\]
\end{theorem}

     \begin{figure}
   \centering
{\setlength\fboxsep{0pt}{\includegraphics[width=.3\columnwidth]{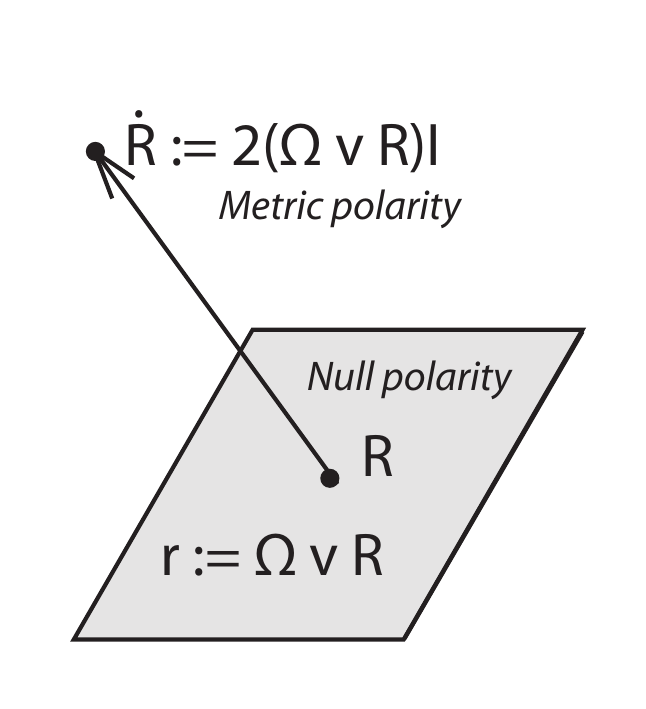}}}
{\setlength\fboxsep{0pt}{\includegraphics[width=.65\columnwidth]{threeNullPlanesWithPerpLabeled-02.pdf}}}
\caption{\emph{Left:} The infinitesimal direction of motion of a point is the composition of the null polarity of the velocity state $\velo$ and the metric polarity.\emph{Right}: A similar picture arises at each point of a null line.}
\label{fig:derivAtPoint}
\end{figure}
 
\myboldhead{Direction of motion vector field determined by $\velo$} 
The global velocity state $\velo$ sets every point in space moving in a particular direction.  We show how to find this direction as an example of how PGA kinematics works. 
The theorem implies that $\dot{\vec{P}}(t) = 2(\velo_s(t) \times \vec{P}_s(t))$, where $\vec{A} \times \vec{B} := \vec{A}\vec{B} - \vec{B}\vec{A}$ is the \emph{commutator} of $\vec{A}$ and $\vec{B}$, sometimes called the \emph{Lie bracket} in this context.  Specializing to $t=0$ and dropping the subscripts since the body and space coordinate systems agree there, we obtain $\dot{\vec{P}}(0) = 2(\velo \times \vec{P})$. Now we can apply one of the many identities of the geometric product, this one involving the product of a 2-vector and a 3-vector: $\velo \times \vec{P} = (\velo \vee \vec{P})\eye$. 

To understand the right-hand side we must first introduce the concept of the \emph{null system}, that plays an essential role in  this approach to rigid body mechanics. In \cite{moebius37}, M\"{o}bius  introduced the concept of a  \emph{null line}. In our terminology, a null line $\sigo$  of a global velocity state $\velo$ is a line such that a force acting along the line $\sigo$ is indifferent to  $\velo$, that is, achieves no virtual work in the presence of $\velo$. M\"{o}bius discovered that, in general, the set of all null lines passing through a point $\vec{P}$ lie in a plane of the point, the \emph{null plane} of $\vec{P}$; and the null lines lying in a plane pass through a point of the plane, the \emph{null point} of the plane.  In PGA, a null line satisfies $\sigo \wedge \velo = 0$, the null plane of $\vec{P}$ is $\velo \vee \vec{P}$,  while the null point $\vec{m}$ is $\velo \wedge \vec{m}$.  Together these two operations are part of the \emph{null polarity} defined on the whole algebra.  It is a polarity since it is grade-reversing and is an involution.

To return to our vector field: the right-hand side $(\velo \vee \vec{P})\eye$ is the product of two polarities. First, $\velo \vee \vec{P}$ is the null polarity associated to $\velo$. This is followed by the metric polarity (multiplication by $\eye$), sending an element to its \quot{orthogonal complement}. This gives the perpendicular direction to the plane $\velo \vee \vec{P}$.  But this is exactly the infinitesimal direction of motion of the point $\vec{P}$.  Consult \Fig{fig:derivAtPoint}.
 
\subsection{Rigid body mechanics}

In the interest of space, we focus here on a qualitative treatment of this topic.  Readers interested in full details should consult \cite{gunnthesis}, Ch. 9, \cite{gunn2011}, or \cite{gunnFull2010}.  

\myboldhead{Statics} We begin by rephrasing standard results of statics in this framework.  A \emph{simple force} is represented by a simple 2-vector $\pip$ where $\pip_e$ is the direction vector of the force and  $\pip_i$, considered as a 3-vector, is the \emph{moment} of the force and describes how the line of force is translated away from the origin.  A resultant of a set of forces $\pip_k$ is then given by the sum: $\pip = \sum_k{\pip_k}$.  The forces are in equilibrium $\iff \pip = 0$.  The resultant is a simple force $\iff \pip \wedge \pip = 0$.  Finally, the resultant is a force couple $\iff \pip = \pip_i$, that is, $\pip$ is ideal. Hence, just as  ideal elements allow PGA kinematics to handle translations and rotations uniformly, they do the same in dynamics for forces and force couples: A force couple is a force carried by an ideal line.

\myboldhead{Newtonian particles} Newtonian particles map in a similarly direct way to the bivector framework.  A particle is considered a point mass with mass $m$, located at position $\vec{R}$ and moving with velocity $\dot{\vec{R}}$. Then one defines:
\begin{compactenum}
\item The \emph{spear} of the particle is $\vec{\Lambda} := \vec{R} \vee \dot{\vec{R}}$.
\item The \emph{momentum state} of the particle is $\momo := m\vec{\Lambda}$.
\item The \emph{velocity state} of the particle is $\pvelo := \vec{\Lambda}\eye$. 
%
\item The \emph{kinetic energy} $E$ of the particle is 
$E := \frac{m}{2} \|\dot{\vec{R}}\|_\infty^{2} = -\frac{m}{2} \vec{\Lambda} \cdot \vec{\Lambda} = -\frac{1}{2}\pvelo \vee \momo$. 
\end{compactenum}
Physics proper begins by observing that in the absence of external forces, $\vec{\Lambda}$, $\momo$, $\pvelo$, and $E$ are conserved quantities.

\myboldhead{Rigid body mechanics} We restrict our further comments on mechanics to a bare minimum to arrive at the Euler equations of motion. Readers interested in the details are referred to \cite{gunn2011} or \cite{gunnFull2010}. 
Mimicking standard approaches of rigid body mechanics one can define a rigid body composed of such newtonian particles.  By summation/integration one can then define analogous quantities for the rigid body, such as momentum and kinetic energy.  The \emph{inertia tensor} can then be introduced, which is a quadratic form $\vec{A}$ on the space of velocity states $\bigwedge^2$ (bi-vectors).  $\vec{A}$ acts on vectors in the body coordinate system.  Then $\vec{A}(\velo, \velo)$ is the kinetic energy of the rigid body under the influence of the global velocity state $\velo$. We write $\vec{A}(\velo)$ for the associated polarizing operator, that produces the momentum state $\momo$ of the rigid body corresponding to this velocity.  Using the fact that the momentum in space, in the absence of external forces, is conserved, one can then derive the Euler equations of motion:
\begin{align*}\label{eqn:eqnmot}
\dot{\vec{g}} &= \vec{g}\velo_c\\
\dot{\velo}_c &= 2\inert^{-1}(\inert(\velo_c) \times \velo_c )
\end{align*}

\myboldhead{Analysis}This is a set of ordinary linear differential equations that display very attractive numerical features.  The solution space of each is a 6D quantity, for a total solution space of 12 dimensions. $\vec{g}$ is an element of the 8D even-subalgebra, while $\vec{\velo_c}$ is a bi-vector itself.  Hence the total representation space is 14 dimensions; there is very little space for \quot{wandering off} the solution.  In fact, we have obtained excellent practical results by simply normalizing $\vec{g}$ to have norm 1, a subspace (the so-called rotor group) of dimension 6 within the even subalgebra. Compare this approach to the standard linear algebra approach, where the use of 4x4 (or 3x4) matrices to represent the Euclidean group introduces many extra dimensions to the representation space, forcing the use of Lagrange multiplies or other methods to constrain the solutions. The situation with the alternative geometric algebra model CGA is similar, see \cite{gunn2017a} \S 7.2.  That is, the PGA approach provides a nearly optimally small representation space for the Euler equations, minimizing numerical problems. Furthermore, from the solution one can directly read off the axis of the resultant screw motion $\vec{g}$ or the global velocity state $\velo_c$.  

\onecolumn
\begin{table}
\begin{center}
\normalsize
\setlength{\extrarowheight}{4pt}
\begin{tabular}{| p{.45\columnwidth} | p{.45\columnwidth} |}
\hline
\textbf{Projective geometric algebra} & \textbf{Vector and linear alg + analytic geom.} \\
\hline
\hline
Unified representation for points, lines, and planes based on a graded exterior algebra; all are \quot{equal citizens} in the algebra.
&
The basic primitives are  points and vectors and all other primitives are built up from these. For example, lines in 3D sometimes  parametric, sometimes w/ Pl\"{u}cker coordinates. \\
\hline
exterior algebra provides robust meet and join operators that deal correctly with parallel entities
&
Meet and join operators only possible when homogeneous coordinates are used, even then tend to be \emph{ad hoc} since points have distinguished role and ideal elements  rarely integrated. \\
\hline
Unified, high-level treatment of euclidean (\quot{finite}) and ideal (\quot{infinite}) elements of all dimensions.
&
Points (euclidean) and vectors (ideal) have their own rules, user must keep track of which is which; no higher-dimensional analogues for lines and planes.\\
\hline
Unified repn. of isometries based on sandwich operators which act uniformly on points, lines, and planes.
&
Matrix representation for isometries has different forms for points, lines, and planes.\\
\hline
Same representation for operator and operand: $\vec{m}$ is the  plane as well as the reflection in the plane.
&
Matrix representation for reflection in $\vec{m}$ is different from the vector representing the plane. \\
\hline
Compact, universal expressive formulas and constructions based on geometric product (see \Tab{tab:pga2d} and \Tab{tab:pga3d}) valid for wide range of argument types and dimensions.
&
Formulas and constructions are \emph{ad hoc}, complicated, many special cases, separate formulas for points/lines/planes.\\
\hline
Faithful representation of Lie algebra and Lie group of $\Euc{n}$ is geometrically intuitive and numerically optimal.
&
Representation of Lie algebra and Lie group of $\Euc{n}$ as subgroups of matrix group $GL(n, \mathbb{R})$ is not geometrically intuitive and is numerically problematical. \\
\hline
\end{tabular}
\myvspace{.1in}
\caption{A comparison of PGA with the standard approach. }
\end{center}
\end{table}
\twocolumn
}

\bibliography{GunnRef}

\begin{thebibliography}{Gun17b}

\bibitem[Bos18]{observable}
Mike Bostock.
\newblock Observable notebooks: a reactive {J}ava{S}cript environment, 2018.
\newblock https://observablehq.com.

\bibitem[Cli73]{clifford73}
William Clifford.
\newblock A preliminary sketch of biquaternions.
\newblock {\em Proc. London Math. Soc.}, 4:381--395, 1873.

\bibitem[Cli78]{clifford78}
William Clifford.
\newblock Applications of {G}rassmann's extensive algebra.
\newblock {\em American Journal of Mathematics}, 1(4):pp. 350--358, 1878.

\bibitem[DFM07]{dfm07}
Leo Dorst, Daniel Fontijne, and Stephen Mann.
\newblock {\em Geometric Algebra for Computer Science}.
\newblock Morgan Kaufmann, San Francisco, 2007.

\bibitem[dK17a]{ganjacs}
Steven de~Keninck.
\newblock Ganja coffeeshop, 2017.
\newblock https://enkimute.github.io/ganja.js/examples.

\bibitem[dK17b]{ganja}
Steven de~Keninck.
\newblock Ganja: Geometric algebra for javascript, 2017.
\newblock https://github.com/enkimute/ganja.js.

\bibitem[Gla90]{gg90}
Andrew~S. Glassner.
\newblock Useful 3d geometry.
\newblock In Andrew~S. Glassner, editor, {\em Graphics Gems}, pages 297--300.
  Academic Press, 1990.

\bibitem[Gra44]{grassmann44}
Hermann Grassmann.
\newblock {\em Ausdehnungslehre}.
\newblock Otto Wigand, Leipzig, 1844.

\bibitem[Gun11a]{gunnthesis}
Charles Gunn.
\newblock {\em Geometry, Kinematics, and Rigid Body Mechanics in Cayley-Klein
  Geometries}.
\newblock PhD thesis, Technical University Berlin, 2011.
\newblock \url{http://opus.kobv.de/tuberlin/volltexte/2011/3322}.

\bibitem[Gun11b]{gunn2011}
Charles Gunn.
\newblock On the homogeneous model of euclidean geometry.
\newblock In Leo Dorst and Joan Lasenby, editors, {\em A Guide to Geometric
  Algebra in Practice}, chapter~15, pages 297--327. Springer, 2011.

\bibitem[Gun17a]{gunn2017b}
Charles Gunn.
\newblock Doing euclidean plane geometry using projective geometric algebra.
\newblock {\em Advances in Applied Clifford Algebras}, 27(2):1203--1232, 2017.

\bibitem[Gun17b]{gunn2017a}
Charles Gunn.
\newblock Geometric algebras for euclidean geometry.
\newblock {\em Advances in Applied Clifford Algebras}, 27(1):185--208, 2017.

\bibitem[Hil13]{hildebrand13}
Dieter Hildebrand.
\newblock {\em Fundamentals of Geometric Algebra Computing}.
\newblock Springer, 2013.

\bibitem[HLR01]{hlr01}
David Hestenes, Hongbo Li, and Alyn Rockwood.
\newblock A unified algebraic approach for classical geometries.
\newblock In Gerald Sommer, editor, {\em Geometric Computing with Clifford
  Algebra}, pages 3--27. Springer, 2001.

\bibitem[HS87]{hessob87}
David Hestenes and Garret Sobczyk.
\newblock {\em Clifford Algebra to Geometric Calculus}.
\newblock Fundamental Theories of Physics. Reidel, Dordrecht, 1987.

\bibitem[Lan71]{lang71}
Serge Lang.
\newblock {\em Algebra}.
\newblock Addison-Wesley, 1971.

\bibitem[LLD11]{lasenby2011}
Anthony Lasenby, Robert Lasenby, and Chris Doran.
\newblock Rigid body dynamics and conformal geometric algebra.
\newblock In Leo Dorst and Joan Lasenby, editors, {\em Guide to Geometric
  Algebra in Practice}, chapter~1, pages 3--25. Springer, 2011.

\bibitem[Sel00]{selig00}
Jon Selig.
\newblock Clifford algebra of points, lines, and planes.
\newblock {\em Robotica}, 18:545--556, 2000.

\bibitem[Stu91]{study91}
Eduard Study.
\newblock Von den {B}ewegungen und {U}mlegungen.
\newblock {\em Mathematische Annalen}, 39:441--566, 1891.

\bibitem[Stu03]{study03}
Eduard Study.
\newblock {\em Geometrie der {D}ynamen}.
\newblock Tuebner, Leibzig, 1903.

\bibitem[Wik]{wikiAD}
Wikipedia.
\newblock \url{https://en.wikipedia.org/wiki/Automatic_differentiation}.

\end{thebibliography}
\bibliographystyle{alpha}

\end{document}